\documentstyle{amsppt}
\voffset-1cm
\nologo
\magnification=\magstep1
\leftheadtext{Kenneth A. Ribet}
\rightheadtext{Galois representations and modular forms}

%
\newcount\sectCount
\edef\nonouterhead{\noexpand\head}
\def\section#1\endhead{\advance\sectCount by 1%
\nonouterhead\the\sectCount. #1\endhead\noindent}


\newcount\refCount
\def\newref#1 {\advance\refCount by 1
\expandafter\edef\csname#1\endcsname{\the\refCount}}


  \newref  ASHGROSS
  \newref  AL   
  \newref  GP  
  \newref  ANTWERPIV 
  \newref  BOSMA 
  \newref  BOSTONSURVEY
  \newref  BOSTONGRANVILLE
  \newref  BUZ
  \newref  CAREARLY
  \newref  COATESIP
  \newref  CONNELL 
  \newref  CORNELLSILVERMAN
  \newref  COX
  \newref  DARMON
  \newref  NEWDARMON
  \newref  DELIGNESERRE
  \newref  DIAMONDREF
  \newref  DIAMONDNEW
  \newref  FALT
  \newref  FLACH
  \newref  FONT 
  \newref  FONTAINEMAZUR
  \newref  FREYONE
  \newref  FREYTWO
  \newref  GELBART
  \newref  LABESSE 
  \newref  GOUCONTROL
  \newref  GOUVEA
  \newref  GOUSURVEY
  \newref  HAYESRIBET
  \newref  HEARSTRIBET
  \newref  HINDRY
  \newref  AXJ 
  \newref  ANOTHERAXJ
  \newref  KNAPP
  \newref  LANGBOOK
  \newref  LANGABELIANVAR
  \newref  LANGFILE
  \newref  LANG  
  \newref  LANGLANDS
  \newref  HWL
  \newref  LI  
  \newref  MATSUMURA
  \newref  SILVERMAN 
  \newref  EIS  
  \newref  PRIMEDEGREE
  \newref  MAZUR  
  \newref  GAD
  \newref  MAZURNOTES
  \newref  MAZURNUMBER
  \newref  MAZURTILOUINE
  \newref  MIYAKEOLD
  \newref  MIYAKE
  \newref  MURTY
  \newref  MURTYASED
  \newref  OESTERLE
  \newref  OGG
  \newref  PRASAD
  \newref  RAMA 
  \newref  RIBETSURVEY
  \newref  ICM
  \newref  FERMAT
  \newref  TOULOUSE
  \newref  KOREA
  \newref  NOTICES
  \newref  VIDEO
  \newref  MOTIVES
  \newref  RUBINSILVER
  \newref  RSNOTE
  \newref  SERRETAU
  \newref  SERREOLD
  \newref  SERREIM
  \newref  SERRECOURS
  \newref  ARCATA
  \newref  DUKE
  \newref  GROUPESALG
  \newref  ST
  \newref  SHIMURAOLD
  \newref  SHCRELLE
  \newref  SHIMURABOOK
  \newref  SHIMURAPINK
  \newref  SHIMURAHECKE
  \newref  SHIMURAFACTORS
  \newref  SILVERMANOLD
  \newref  SILVERMANNEW
  \newref  TATE  
  \newref  TATELETTER
  \newref  TW  
  \newref  TUNNELL
  \newref  MvS
  \newref  WEIL
  \newref  WILES


\let\phi\varphi
\def\U#1{\Z/#1\Z} 
\def\UU#1{(\Z/#1\Z)^\ast} 
\let\a\alpha
\def\Q{{\bold Q}}
\def\Z{{\bold Z}}
\def\C{{\bold C}}
\def\Qbar{\overline{\Q}}
\def\Qp{\Q_p}
\def\Qpbar{\Qbar_p}
\def\F{{\bold F}}

\def\Fp{\F_p}

\def\T{{\Bbb T}}

\def\xo#1{X_0(#1)}
\def\go#1{\Gamma_0(#1)}
\def\ssigma{{}^\sigma\mkern-\thinmuskip}
\def\hlinefill{\leaders\hrule height 3pt depth -2.5pt\hfill}
\def\emrule{\thinspace\hbox to 0.75em{\hlinefill}\thinspace}
\def\GalQ{\Gal(\Qbar/\Q)}
\def\seteq{\mathrel{:=}}
\def\SL#1{{\bold{SL}}(2,#1)}
\def\GL#1{{\bold{GL}}(2,#1)}
\def\GLTWO{{\bold{GL}}(2)}
\def\varGLTWO{\bold{GL}_2}

\def\Weil{Taniyama-Shimura}

\def\notdivide{\setbox1=\hbox{$|$\llap{
\hbox{/}\kern-0.75pt}}\mathrel{\box1}}

\def\<<{\leavevmode
  \raise0.28ex\hbox{$\scriptscriptstyle\langle\!\langle$}\nobreak
  \hskip -.6pt plus.3pt minus.2pt$\,$}
\def\>>{$\,$\nobreak\hskip -.6pt plus.3pt minus.2pt
  \raise0.28ex\hbox{$\scriptscriptstyle\rangle\!\rangle$}}

\let\tempcirc=\circ
\def\circ{\mathord{\raise0.25ex\hbox{$\scriptscriptstyle\tempcirc$}}}

\def\newop#1
{\expandafter\def\csname #1\endcsname{\mathop{\roman{#1}}\nolimits}}
\newop Aut
\newop Gal 
\newop End 
\newop Frob
\newop tr
\newop Hom

\topmatter

\title Galois representations and modular forms\endtitle
\author Kenneth A. Ribet\endauthor
\affil University of California, Berkeley\endaffil
\address{UC Mathematics Department,
  Berkeley, CA 94720-3840 USA}\endaddress
\email ribet\@math.berkeley.edu \endemail
\abstract\nofrills
In this article, I discuss material which is
related to the recent proof of Fermat's Last Theorem: elliptic curves,
modular forms, Galois representations
and their deformations, Frey's construction, and the
conjectures of Serre and of \Weil.
\endabstract
\thanks 
This article was prepared while I was a research professor at MSRI, where
research is supported in part by NSF grant DMS-9022140.
My work was
to some extent
supported by NSF Grant DMS 93-06898.
I wish to thank
N.~Boston,
H.~Darmon,
F.~Diamond,
F.\,Q.~Gouv\^ea,
B.~Mazur,
V.\,K.~Murty,
C.~O'Neil,
S.~Ribet,  
A.~Silverberg,
R.~Taylor
and A.~Wilkinson
for feedback on drafts of these
notes.
I am especially grateful to B.~Conrad and L.~Goldberg
for
detailed comments.
\endthanks

\endtopmatter

\document
\hrule height0pt \vskip-16pt \hrule height0pt
\section Introduction\endhead
This article is a revised version of the notes which were
distributed at my Progress in Mathematics
lecture at the August, 1994 Minneapolis Mathfest.  When I was
first
approached in 1993 by the Progress in Mathematics committee,
I was
asked to discuss
``the mathematics behind Andrew Wiles's solution of the
Fermat conjecture.''
As the reader is no doubt aware, 
Wiles had announced 
in a series of
lectures at the Isaac Newton Institute for Mathematical Sciences
that he was able to prove 
for semistable elliptic curves the conjecture
of Shimura and Taniyama to the effect that all elliptic
curves over~$\Q$ are ``modular'' (i.e., attached
to modular forms in a sense which will be explained below).
Fermat's Last Theorem follows from this result,
together with a theorem that I proved
seven years
ago~\cite\FERMAT.  

At the time of the Mathfest, however, a gap which
had appeared in
Wiles's work had not yet been repaired
(see \cite\AXJ)\emrule
Fermat's ``Theorem'' was still a conjecture.
Nevertheless, it was readily apparent
that the methods introduced
by Wiles were significant and deserving of
attention.  Most notably, these
methods had been used to construct for the first
time an infinite set $\Cal E$ of modular elliptic curves
over~$\Q$ with the following property:
if $E_1$ and~$E_2$ are elements of~$\Cal E$, then
$E_1$ and $E_2$ are non-isomorphic even when viewed
as elliptic curves over the complex field~$\C$
(see~\cite\RSNOTE).  In my lecture at the Mathfest,
I stressed this achievement of Wiles and discussed the
analogy between the \Weil\ conjecture and 
conjectures of~Serre~\cite\DUKE\ about two-dimensional
Galois representations.

As I began to
revise these notes, I found that the situation had changed
dramatically for the better.
Wiles announced
in October, 1994
that the bound he sought in
his original proof
could be
obtained by a method which circumvented the Euler system
construction~\cite\ANOTHERAXJ.  
This method arose from an observation that Wiles made
early in his investigation:
the 
required
upper bound would follow from a proof that certain
Hecke algebras are complete intersection rings.  
This statement about Hecke algebras
is the main theorem of
an article written jointly by
Richard Taylor and Andrew Wiles~\cite\TW;
these authors announced their result at the same time that
Wiles disseminated a revised version of his original
manuscript~\cite\WILES.

It is premature to undertake a detailed analysis
of the results of Wiles and Taylor-Wiles.
The aim of this survey is more modest:
to present an introduction to the circle of
ideas which form the background for these results.
Because of the intense publicity surrounding Fermat's Last
Theorem, a good deal of the material I have chosen
has been discussed in news and
expository articles which were
written in connection with Wiles's 1993 announcement.
Among these are the author's news item in
the Notices~\cite\NOTICES\ and his article
with Brian Hayes in American Scientist~\cite\HAYESRIBET,
two pieces in the American Mathematical Monthly
\cite{\COX, \GOUVEA}, the report by K.~Rubin
and A.~Silverberg in this Bulletin~\cite\RUBINSILVER,
a long survey by H.~Darmon~\cite\DARMON,
an elementary formulation of the Langlands program
by A.~Ash and R.~Gross~\cite\ASHGROSS,
and an article for undergraduates and
their teachers by N.~Boston~\cite\BOSTONSURVEY.
Books containing articles related to Fermat's Last
Theorem are soon to appear~\cite{\COATESIP, \MURTYASED}.
Furthermore, two videotapes related to Fermat have
been in circulation for some time:
a 1993 lecture by the author is available from
the AMS~\cite\VIDEO, and the Mathematical Sciences
Research Institute has been distributing a videotape based
on the July, 1993 ``Fermat Fest'' which was organized
by the MSRI and the San Francisco Exploratorium.
Readers may also consult material
available by gopher from {\tt e-math.ams.org}.

In view of the burgeoning literature in this
subject, I imagined
these notes mostly as a somewhat biased guide
to reference works and
expository articles.
In the end, what I have written might best be
characterized as an
abbreviated survey with a disproportionately large
list of~references.

\section Modular forms\endhead
We begin by summarizing some background 
material concerning
modular forms.  
Among reference books in the subject, one might
cite \cite\LANGBOOK, \cite\MIYAKE, and~\cite\SHIMURABOOK.
Also,
several books on elliptic curves contain substantial 
material on modular forms; in particular, Knapp's book~\cite\KNAPP\
has been recommended to the author
as a good source for an overview of the
Eichler-Shimura theory relating elliptic curves to certain
modular forms.
For a first introduction to modular forms, a fine
starting point is~\cite{\SERRECOURS, Ch.~VII}.

The modular forms to be considered are 
``cusp forms of weight two on~$\go N$,''
for some integer $N\ge1$.
Here, $\go N$ is
the group of integer matrices
with determinant~1 which are upper-triangular mod~$N$.
A cusp form on this group is, in particular, a
holomorphic
function
on the upper half-plane $\Cal H$ consisting
of complex numbers
with positive imaginary part.  

These functions
are usually presented as converging Fourier series
$\sum_{n=1}^\infty a_n q^n$, where $q\seteq e^{2\pi i z}$.
For the forms
which most interest us,
the complex numbers $a_n$
are algebraic
integers; frequently they are even ordinary integers.

The weight-two cusp forms on~$\go N$
are holomorphic functions $f(z)$
on the complex upper half-plane $\Cal H$.
On requires principally
that $f(z)\,dz$ be invariant under~$\go N$,
i.e., that
$f(z)$ satisfy
the functional equation
$$f\left({az+b\over cz+d}\right)= (cz+d)^2 f(z)$$
for all
$\pmatrix a & b\cr c& d\cr\endpmatrix\in\go N$.
(In particular, one has $f(z+b)=f(z)$ 
for all integers~$b$.)
In addition to the holomorphy and the functional
equation, one imposes
subsidiary conditions at infinity~\cite{\SHIMURABOOK, Ch.~2}.

The group $\go N$
acts on~$\Cal H$ by fractional linear
transformations, with
$\pmatrix a & b\cr
c& d\cr\endpmatrix$ acting as $$z\mapsto {az+b\over cz+d}.$$
The invariance of~$f(z)\,dz$ means that $f(z)\,dz$
arises by pullback from a differential on
the quotient
$Y_0(N)\seteq \go N\backslash\Cal H$.
This quotient is 
a non-compact Riemann surface with a standard compactification,
known as~$\xo N$.  The complement of~$Y_0(N)$ in~$\xo N$
is a finite set, the set of
{\it cusps} of~$\xo N$.
The conditions on infinity satisfied by a cusp form require
the differential on~$Y_0(N)$ associated to~$f(z)$ to extend
to~$\xo N$.  This means simply that
 $f(z)$ is required to vanish at
the cusps.

Let us identify the space $S(N)$ of weight-two
cusp forms on~$\go N$ 
with
the space of
holomorphic differentials on the Riemann surface~$\xo N$.
Then  $S(N)$ has finite dimension:
$\dim S(N)$ is the genus of the curve~$\xo N$.
This integer
may be calculated easily via the Riemann-Hurwitz
formula (applied to the covering $\xo N\to \xo 1$
which corresponds to the inclusion of $\go N$ in~$\go1$).
For example, $\dim S(N)$
is~zero for $N\le10$ and is~one for~$N=11$.
(See \cite{\SHIMURABOOK, Ch.~2}.)

In order to present an example of
a non-zero cusp form, we will exhibit a
non-zero
element of the 1-dimensional space
$S(11)$.
Namely,
consider
the formal power series 
with integral coefficients
$\sum a_n X^n$
which is defined by the identity
$$ X \prod_{m=1}^\infty (1-X^m)^2(1-X^{11m})^2 =
\sum_{n=1}^\infty a_nX^n.$$
It can be shown that the
holomorphic function $\sum a_n q^n $
(with $q=e^{2\pi i z}$)
is a cusp form $h$
of weight two on~$\go{11}$, i.e., a generator of
the 1-dimensional space~$S(11)$, cf.~\cite{\SHIMURABOOK,
Example~2.28}.
(One has
$h(z)=\eta(z)^2\eta(11z)^2$, where $\eta(z)$, the Dedekind
$\eta$-function,
is the standard example of a modular form of weight~$1/2$.)

For each integer $n\ge1$, the $n$th {\it Hecke operator\/}
on~$S(N)$ is an endomorphism
$T_n$ of~$S(N)$, whose action is generally written on the right:
$f\mapsto f|T_n$. 
The various $T_n$ commute with each other and are interrelated
by identities which express a given $T_n$ in terms of the
Hecke operators indexed by the prime factors of~$n$.
If $p$ is a prime, the operator $T_p$ 
is defined as a composite $\beta\circ\alpha$, where 
$\alpha\:S(N)\to S(Np)$ and $\beta\:S(Np)\to S(N)$ are
standard homomorphisms, known as ``degeneracy maps.''
In terms of Fourier coefficients,
$T_p$ is given by the following rules:
If $p$ is not a divisor of~$N$ and $f=\sum a_n q^n$,
then $f|T_p$ has the Fourier expansion 
$$ \sum_{n=1}^\infty a_{np} q^n + p \sum_{n=1}^\infty a_nq^{pn}.$$
If $p$ divides~$N$, then
$f|T_p$ is given by
$\sum_{n=1}^\infty a_{np} q^n$.  
The $T_n$ have a strong tendency to be diagonalizable on~$S(N)$:
for $n$ prime to~$N$, $T_n$ is self-adjoint with respect to the
Petersson pairing on~$S(N)$, and is therefore semisimple
\cite{\SHIMURABOOK, \S3.5}.

The elements of~$S(N)$ having special
arithmetic interest are the {\it normalized
eigenforms\/} in~$S(N)$; these are the non-zero
cusp forms $f=\sum a_nq^n$
which are eigenvectors for all the $T_n$ and which
satisfy the normalizing condition $a_1=1$.
(The latter condition is imposed mostly for convenience,
since an arbitrary eigenform is a multiple of a 
normalized one.)
If $f$ is such an eigenform, its Fourier coefficients
and its eigenvalues coincide: one has
$f|T_n = a_n f$ for all $n\ge1$.  The
normalized eigenforms are
``arithmetic'' because the~$a_n$ belong to the realm
of algebraic number theory.  In fact,
if $f$ is a normalized
eigenform, then the subfield of~$\C$
generated by the $a_n$ is a finite (algebraic) extension $E$
of~$\Q$ in~$\C$, and the elements $a_n$ of~$E$ are
algebraic integers (see
\cite{\SHIMURABOOK, Th.~3.52} or~\cite{\DELIGNESERRE,
Prop.~2.7}).

It is a mildly complicating
fact that the normalized eigenforms in~$S(N)$
do not always form a basis of~$S(N)$.  In other words,
the commuting
operators $T_n$ are not necessarily 
all diagonalizable.
This problem, which arises from those $n$ which have a
common factor with~$N$, can be repaired by the
introduction of {\it newforms\/}~\cite\AL.  Briefly, a newform
is a normalized eigenform $f=\sum a_n q^n$ for which
the space
$$\Big\{\, g\in S(N) \, \Big| \, g|T_n = a_n g \hbox{ for
all $n$ prime to $N$}\,\Big\}$$
contains only $f$ and its multiples.
Atkin and Lehner showed in~1970
that $S(N)$ has a basis built out of 
suitable transforms of
the newforms
in the spaces $S(M)$,
where $M$ runs over the positive divisors of~$N$. 
Namely,
for each~$M$,
there are natural embeddings of the space $S(M)$ into~$S(N)$
which are indexed by the positive divisors of~$N/M$.
Namely,  
to each divisor $d$ of $N/M$ is associated the embedding
$\delta_d\:S(M)\hookrightarrow S(N)$ which sends 
$\sum a_nq^n$ to the series $\sum a_nq^{nd}$.
The
embeddings $\delta_d$ are {\it degeneracy mappings}; one
such mapping was introduced briefly above in connection with
the definition of the $p$th Hecke operator~$T_p$.  With the help
of these
degeneracy mappings, the full space
$S(N)$
may be reconstructed purely in terms of the newforms
in the spaces $S(M)$ with $M$ dividing~$N$.
Specifically, let $S(M)^{\roman{new}}$ be the subspace of~$S(M)$
spanned by its newforms.  Then 
$$  S(N) = \bigoplus_{M|N} \bigoplus_{d|{N\over M}} 
\delta_d\big(S(M)^{\roman{new}}\big)
.$$
We conclude this discussion with the remark that
the work of Atkin and Lehner was
generalized by T.~Miyake~\cite\MIYAKEOLD\ in the setting of
automorphic forms on~$\GLTWO$ and then by W.~Li \cite\LI\
in the setting of classical modular forms.

\section Elliptic curves\endhead
We next introduce some
foundational concepts pertaining to
elliptic curves.  For an in-depth treatment of these concepts,
the reader may consult the large number
of textbooks and
monographs which focus on
elliptic curves.  (As was indicated above, some of these
books
discuss modular forms as well.)
Rather than list these references here,
we refer the reader to the bibliography of a
recent
book review written by W.~R.~Hearst~III
and the author~\cite\HEARSTRIBET.
One book which has appeared since that review was written
is the recent
``Advanced Topics in the Arithmetic of Elliptic Curves''
by J.~H.~Silverman~\cite\SILVERMANNEW.

Elliptic curves are distinguished by the fact that they are
simultaneously {\it curves}, i.e., varieties of dimension~1,
and {\it abelian varieties}, projective varieties which
are endowed with a group law.  The definition is very simple:
an elliptic curve over a field~$k$ is a projective
non-singular curve $A$ of genus~1 over~$k$ which is furnished
with a distinguished rational point~$O$.  An elliptic curve
is thus a pair $(A,O)$; an isomorphism between
two pairs $(A,O)$ and $(A',O')$
is an isomorphism $A\buildrel\sim\over\to A'$ which maps $O$
to~$O'$.  By convention, one usually omits explicit
mention of the point $O$ and refers simply to an elliptic
curve $A$ over~$k$.

One obtains elliptic curves from 5-tuples 
of elements of~$k$: to $(\a_1, \a_2, \a_3, \a_4, \a_6)$
one associates the projective curve whose affine equation
is the generalized Weierstra\ss\ equation
$$y^2 + \a_1xy + \a_3 y = x^3 + \a_2 x^2 + \a_4 x + \a_6.$$
This curve has genus~1 if it is non-singular;
the non-singularity is detected by the non-vanishing of
the {\it discriminant\/} of the equation.  (The discriminant
is a polynomial in the $\a_i$ which
will not be reproduced here; it can be
found, for instance, on page~36 of~\cite\TATE.)
A short calculation shows that the plane curve with the
indicated equation has exactly one point which does not
lie on the affine plane; we take this point to be the distinguished
point~$O$.  Using the Riemann-Roch theorem, one checks
that every elliptic curve over~$k$ is isomorphic
to one obtained by this process. (See, e.g., 
\cite{\GROUPESALG, Ch.~II} for a discussion of the Riemann-Roch theorem
and
\cite{\SILVERMANOLD,
Ch.~III, \S3} for a proof that elliptic curves are
given by generalized Weierstra\ss\ equations.)

As we have suggested, an elliptic curve
$A$ over~$k$
may be viewed as
a commutative {\it algebraic group\/} over~$k$.
Concretely, suppose that one has chosen
a Weierstra\ss\ equation for~$A$ as above.
For each field
$K$ containing~$k$,
let $A(K)$ be the
subset of the projective plane
over~$K$ which is defined by the Weierstra\ss\
equation.
This set is endowed with
a classical group law, often known as the
``chord and tangent operation.''  
In this law, $O$ is the identity element of~$A(K)$,
and three distinct elements of~$A(K)$ sum to~$O$ if and only
if they are collinear.
The composition law on~$A(K)$ can
be described explicitly
in terms of coordinates by a family of polynomial
equations with coefficients in~$k$;
this family depends on the chosen Weierstra\ss\
equation, but is independent of~$K$.
(For a recent discussion concerning families of such
equations, see~\cite\BOSMA.)

If $A$ is an elliptic curve over~$\Q$,
then one can choose the Weierstra\ss\
coefficients $\a_i$
to be rational {\it integers}.
The coefficients are essentially unique if we demand
that the discriminant of the equation have the smallest
possible value;
the discriminant
is then said to be the
(minimal) discriminant $\Delta$ of~$A$.

Those prime numbers which divide $\Delta$
are the primes at which $A$ has {\it bad reduction}; the
others are the primes at which $A$ has {\it good
reduction}.  The point is that,
for good primes~$p$,
the minimal equation, when viewed
mod~$p$, yields an elliptic curve $\tilde{A}_p$
over the
finite field $\Z/p\Z$.
The reduced curve
$\tilde A_p$ over~$\Z/p\Z$
is a projective plane curve over~$\Z/p\Z$; 
the group $\tilde{A}_p(\Z/p\Z)$ 
consisting of the points of~$\tilde A_p$ with coordinates
in~$\Z/p\Z$
is then a finite
abelian group.
A theorem of H.~Hasse states that the 
order of this group is approximately~$p$.
More precisely, the
integer $$b_p\seteq
p{+}1-\#(\tilde{A}_p(\Z/p\Z))$$
is bounded in absolute value by~$2\sqrt p$
 \cite{\SILVERMANOLD,
Ch.~V, Th.~1.1}.

To work with a concrete example, let us introduce the
elliptic curve $A$ over~$\Q$
defined by the equation $y^2+y=x^3-x^2$.
The group $A(\Q)$ has five evident points: the origin~$O$,
and the four affine
points gotten by taking $y\in\{0, -1\}$ and $x\in\{0,1\}$.  One may
verify
that these points form a {\it subgroup\/}
of~$A(\Q)$.  Since they remain distinct under
reduction mod~$p$, $\tilde{A}_p(\Z/p\Z)$ has a subgroup of
order~5 whenever $A$ has good reduction at~$p$.  Accordingly,
the number $\#(\tilde{A}_p(\Z/p\Z))$ is divisible by~5, so
that $b_p\equiv p{+}1$ mod~5 for all primes $p$
which do not divide the discriminant of~$A$.

Now Table~1 of~\cite\ANTWERPIV\ informs us that
the minimal
discriminant of~$A$ is $-11$.  Therefore, integers $b_p$
are defined for all $p\ne11$.  When $p$ is small, it is
not hard to compute~$b_p$; one finds
$b_2=-2$, $b_3=-1$, $b_5=1$, and so on.  
(The mod~5 congruence on the~$b_p$,
plus Hasse's bound $|b_p|<2\sqrt p$, determines
the first few of these integers.)
Anticipating the
conjecture 
of Shimura and Taniyama
which will be discussed below,
I would like to point out
that the arithmetically defined
$b_p$ are known to agree with the 
prime-indexed
coefficients $a_p$ of the weight-two modular form
$h$
on~$\go{11}$ which was introduced above.
The astounding identity $b_p=a_p$
is a special case of a
relation which was
discovered by Eicher and Shimura (and which is 
recalled in~\cite\SHCRELLE).
A priori, the Eichler-Shimura
relation gives information only 
about those elliptic
curves which 
are related geometrically to modular curves of the
form $\xo N$.
The \Weil\ conjecture states that
{\it all\/}
elliptic curves over~$\Q$ can be so related.

Next, suppose that $A$ and~$A'$ are elliptic curves
over~$k$.
An {\it isogeny\/} $A\to A'$ over~$k$
is a non-constant map of curves over~$k$
which takes the distinguished point $O$ of~$A$ to the
corresponding point $O'$ of~$A$.
If such a map
exists, 
then it is a map of algebraic groups, and
$A$ and~$A'$ are said to be isogenous
over~$k$.  (As the terminology suggests, the relation
of being isogenous is symmetric; it is in fact an
equivalence relation.)  It is not difficult to show
that two isogenous elliptic curves over~$\Q$ have
the same primes of bad reduction.  (This follows from
the criterion of N\'eron-Ogg-Shafarevich\emrule
see \cite\ST\ or \cite{\SILVERMANOLD, Ch.~VII, \S7}.)
Moreover, if $A$ and~$A'$
are isogenous over~$\Q$, then for each good reduction
prime~$p$, one has  
$$\#(\tilde{A}_p(\Z/p\Z))=\#(\tilde{A'}_p(\Z/p\Z)).$$
In other words, the numbers $b_p$ are the same for $A$
and~$A'$, so that one might make the informal
statement
that $A$ and~$A'$ are ``equivalent
arithmetically.''
A striking theorem of G.~Faltings~\cite{\FALT,\S5} states
that, conversely, two elliptic curves $A$ and~$A'$
over~$\Q$ are isogenous if
$\#(\tilde{A}_p(\Z/p\Z))=\#(\tilde{A'}_p(\Z/p\Z))$
for all primes $p$ at which the two curves have
good reduction.  
(In~\cite\FALT, Faltings proves a more general statement
about homomorphisms between abelian varieties over
number fields, thus confirming conjectures of J.~Tate.) 

Each elliptic curve $A$ over~$\Q$ has a
{\it conductor}, which is a positive
integer divisible precisely by the
primes at which $A$ has bad reduction.
For example, the curve with equation $y^2+y=x^3-x^2$
has conductor~11.
Although the definition of the conductor is somewhat
unenlightening
(see for example~\cite{\SERREOLD, \S2}),
a well known algorithm of J.~Tate~\cite\TATE\
makes it possible to compute conductors
by hand in specific cases.  
Alternatively,
the 
computer algebra program {\tt gp}~\cite\GP,
the
Mathematica package~\cite\SILVERMAN\ and the
Maple package {\tt apecs} all
implement Tate's algorithm\emrule and more generally make
it easy to perform elliptic curve calculations 
on a workstation or personal computer.
(The latter package is recommended by F.~Gouv\^ea.)

Because the conductor is divisible exactly by
the set of ``bad primes,'' the
conductor and the minimal discriminant of~$A$
are divisible by the same set of
prime numbers.  
Nevertheless, these integers have a completely
different feel.
For one thing, the conductor of an elliptic 
curve
depends only on its $\Q$-isogeny class, while the
discriminant may change under isogeny.
(The curve with equation $y^2+y=x^3-x^2-10x-20$
is isogenous to the conductor-11 curve that we
have been discussing; its conductor is~11, but its
minimal discriminant is~$-11^5$.)
For another, the conductor is always positive;
the minimal discriminant
may be positive or negative.  Finally,
the latter number may be
divisible by a large power of a given prime, whereas
the conductor of an elliptic curve tends to be divisible
only by low powers of its prime divisors. 
(In particular,
a prime $p\ge5$ can divide the conductor
at most to the second
power.)
A beautiful formula of A.~Ogg~\cite\OGG\
expresses the conductor
of a given curve $A$ in terms of the minimal
discriminant and the N\'eron
model of~$A$.  (The latter is a curve over~$\Z$ which
may be regarded as the ``best possible'' model for~$A$.
See also \cite{\SILVERMANOLD, App.~C, \S16} for a statement
of Ogg's formula.)

Suppose that the conductor of~$A$ is~$N$.  Then $A$ is said
to be {\it semistable\/} at a prime number~$p$ if $p^2$ does
not divide~$N$.  This means either that $p$ is prime to~$N$,
in which case $A$ has good reduction at~$p$, or else that
$p$ ``exactly'' divides~$N$, in which case the reduction of~$A$
at~$p$ is bad\emrule but not too bad (it is said to be
multiplicative).  If $A$ is semistable at all primes~$p$, i.e., if
$N$ is square free, then the elliptic curve $A$ is said to be
semistable.  Semistable elliptic curves occur in 
connection with Fermat's Last Theorem, and in
other applications.

This is a good point to insert a short digression about
abelian varieties, whose arithmetic
theory is exposed in various
chapters of~\cite\CORNELLSILVERMAN\ and
the recent graduate-level text~\cite\MURTY; see also
the books by such authors as
A.~Weil, S.~Lang~\cite\LANGABELIANVAR, D.~Mumford,
H.~P.~F.~Swinnerton-Dyer, 
G. Kempf
and H.~Lange.
As was mentioned very briefly above, an
abelian variety over a field~$k$ is a
projective algebraic variety
over~$k$ which is equipped with the structure of
an algebraic group over~$k$.  
The first examples of abelian varieties are obtained
by taking Jacobians of curves; the Jacobian of a curve
of genus~$g$ over~$k$ is an abelian variety over~$k$
of dimension~$g$.  
(For the construction of Jacobians, the reader may consult,
e.g., Chapter~V of~\cite\GROUPESALG.)
If $A$ is an elliptic curve over~$k$,
then the chosen origin $O$ of~$A$ enables one to identify
$A$ (endowed with its group structure) with the Jacobian
of~$A$.  In fact, an elliptic curve over~$k$
is nothing other than an abelian variety over~$k$ of dimension~one.
Since an abelian variety of dimension~one is an
elliptic curve (and since
an abelian variety of dimension~zero is just
a single point), one might describe abelian varieties as
``higher-dimensional
analogues'' of elliptic curves.

Abelian varieties play an inevitable role in
our story because the conjecture
which we are about to discuss (a priori one
involving
elliptic curves over~$\Q$ and weight-two newforms
with integer coefficients)
extends naturally to a conjectural dictionary between
arbitrary newforms of weight~two and a certain class 
of abelian varieties over~$\Q$ (see Conjecture~2 below).

\section The \Weil\ conjecture\endhead
The conjecture of Shimura and Taniyama
relates elliptic curves
over~$\Q$ and certain modular forms.
Its history is the subject 
of an engrossing ``file''
compiled by S.~Lang~\cite\LANGFILE.
If I understand correctly, the conjecture was first
posed as a tentative question by Taniyama at the
Tokyo-Nikko conference of~1955.  Shimura stated the
conjecture in its present form in the early 1960s.
In 1968, A.~Weil proved a theorem to the effect
that certain necessary conditions for
an elliptic curve to be ``modular'' are, in fact,
sufficient~\cite\WEIL.  Because of this theorem,
Weil's name has sometimes been attached to the conjecture.

Before giving a formal description of the conjecture,
we refer once again to the elliptic curve $y^2+y=x^3-x^2$,
which has
conductor~11.  The integers $b_p$ which control the numbers
of mod~$p$ points on this curve might be viewed initially
as mysterious quantities for which we seek a ``formula.''
The relation $b_p=a_p$, where $a_p$ is the $p$th coefficient
of the normalized eigenform in~$S(11)$ may be regarded
as such a formula.

The \Weil\ conjecture affirms
that there is an analogous relation for {\it all\/} elliptic
curves over~$\Q$.  Namely, if $A$ is an elliptic curve over~$\Q$
of conductor~$N$, then
one conjectures that
there is a newform $f=\sum a_nq^n$
in~$S(N)$ such that $b_p=a_p$ for all primes $p\notdivide N$.
All coefficients $a_n$ of~$f$ are then 
necessarily
rational integers.
(For another
account of the conjecture, see the article
by K.~Rubin and A.~Silverberg~\cite{\RUBINSILVER, \S1}.)

A geometric formulation of the \Weil\ conjecture may
be given in terms of
the construction presented
in Chapter~7 of~\cite\SHIMURABOOK\
and, from a different perspective, in~\cite\SHIMURAFACTORS.
Suppose 
that $f\in S(N)$ is a normalized eigenform and let
$E$ be the field generated by the coefficients of~$f$.
Shimura associates to~$f$ an abelian variety $A_f$
over~$\Q$
whose dimension is the degree $[E\:\Q]$ and whose arithmetic
incorporates the eigenvalues $a_p$ for $p$ prime to~$N$.
Although the construction of~$A_f$ is perfectly precise,
it is fruitful in this context to regard the
association $f\mapstochar\rightsquigarrow A_f$ as a flabby one
linking $f$ to a ``clump'' of isogenous
abelian varieties, rather than a specific abelian variety
which is singled out up to isomorphism. 
(The relevant notion of isogeny is an appropriate generalization
of the notion of isogeny for elliptic curves.)

If the Fourier coefficients of~$f$ are rational integers,
then $E=\Q$ and $A_f$ has dimension~1.  This means that $A_f$
is an elliptic curve over~$\Q$, which we shall regard
as being defined only up to isogeny. 
According to a theorem of Eichler and Shimura, the
eigenvalues $a_p$ are reflected in the arithmetic
of the elliptic curve $A=A_f$ in the following way.
If $p$ does not divide~$N$, then $A$ has good reduction
at~$p$.  Further, for such $p$,
the $p$th
Fourier coefficient of~$f$ coincides with the quantity
$b_p\seteq p{+}1-\#(\tilde{A}_p(\Z/p\Z))$.
In other words, one has
$b_p=a_p$ for all $p\notdivide N$.
(If $A_f$ has dimension greater than~1, the relation
between $f$ and $A_f$ is a bit more complicated to
formulate, but involves no new ideas.)

In its geometric form, the
\Weil\ conjecture states that
every elliptic curve
$A$ over~$\Q$
is {\it modular\/} in the sense that
it is isogenous to
some curve~$A_f$.
In other words, the conjecture asserts
the surjectivity of
the construction
$f\mapstochar\rightsquigarrow A_f$,
viewed as a map from eigenforms with integral
coefficients (in some~$S(N)$) to isogeny classes of elliptic
curves over~$\Q$.  In analogy with
the arithmetic formulation, when $A$
has conductor~$N$, $A$ is conjectured to be isogenous to an~$A_f$
with $f\in S(N)$.

The connection between the two formulations is as follows.
Suppose that
$A$ is an elliptic curve over~$\Q$ of conductor~$N$, and
assume that $A$ is related arithmetically to an eigenform
in~$S(N)$.  Specifically, assume that $f=\sum a_nq^n$ is a
normalized eigenform in~$S(N)$ with integral coefficients
and that
$a_p=b_p$
for all $p$ prime to~$N$.  (The $b_p$ have their usual meaning.)
Let $A'$ be the elliptic curve $A_f$,
and write $b'_p$ for the analogues of the $b_p$ for~$A'$.
Then one has $b'_p = a_p$ for all $p\notdivide N$ by the formula
of Eichler-Shimura, and hence $b_p=B'_p$ for all $p\notdivide N$.
A theorem of Faltings which was quoted above then ensures that
$A$ and~$A'=A_f$ are isogenous over~$\Q$, so that the
geometric form of the \Weil\ conjecture is true for~$A$.
Conversely, if $A$ is isogenous to an~$A_f$, then numbers
$b_p$ computed for~$A$ coincide with their analogues for~$A_f$.
By the Eicher-Shimura formula, these latter numbers are the
prime-indexed coefficients of~$f$.

In connection with the construction $f\rightsquigarrow A_f$,
let us consider the situation where
$f$ is a
normalized eigenform in~$S(N)$ 
with integral Fourier coefficients
but where $f$ is not necessarily a newform.  What is the
conductor of the
elliptic curve~$A=A_f$?
The answer to this question begins with the fact that 
$f$ is necessarily built from a newform in~$S(M)$, for
some unique divisor $M$ of~$N$.  It is then true that
the conductor of~$A$
is precisely this divisor $M$.  This theorem was proved
by H.~Carayol in~\cite\CAREARLY, following work of
Shimura, Igusa, Deligne and Langlands.

A formulation 
of the \Weil\ conjecture
with a completely different flavor is provided
by
Mazur's article~\cite\GAD.
In this article, Mazur rephrases
the conjecture as a statement
about the Riemann surface associated with an elliptic curve
over~$\Q$.
If $A$ is such an elliptic curve, we write $A(\C)$ for this
Riemann surface, which may be realized as the subset of the complex
projective plane which is defined by a Weierstra\ss\ equation
for~$A$.  This surface is
holomorphically a complex torus.
For each integer~$N$,
consider the subgroup $\Gamma_1(N)$ 
of~$\go N$
consisting of
matrices 
$\pmatrix a & b\cr c&d\endpmatrix\in\go N$ for which
$a\equiv d\equiv 1$~mod~$N$.
By considering $\Gamma_1(N)$ in place of~$\go N$, one
obtains an analogue of~$\xo N$ which is called~$X_1(N)$.
Mazur shows that an elliptic curve $A$ over~$\Q$ is modular
if and only if there exists a non-constant holomorphic
map from $X_1(N)$ to~$A(\C)$, for some positive integer~$N$.
As Mazur explains in his article, one can paraphrase this
condition as the statement that the arithmetic object $A$
is {\it hyperbolically uniformized}, since $\Cal H$ is 
a model for the hyperbolic plane.
This circumstance has led to the charge that Wiles's
proof of Fermat's Last Theorem could not possibly be correct
since its logical structure involves a statement
that may be interpreted in terms of hyperbolic geometry~\cite\MvS.
For a refutation, see~\cite\BOSTONGRANVILLE.

Mazur's observation led Serre to ask
for a
description of
the set of elliptic curves over
the complex field which 
satisfy Mazur's condition:
Which elliptic curves $A$ over~$\C$ 
are modular in the sense that one
can find
a non-constant holomorphic
map $X_1(N)\to A(\C)$ for some positive integer~$N$?
(According to Mazur's theorem, an elliptic curve over~$\Q$
is modular in the usual sense if and only if it is
modular in this new sense when viewed over~$\C$.)
In~\cite\KOREA, I 
provide a conjectural
answer to Serre's question.
Namely, I show that the conjectures made by Serre in~\cite\DUKE\
imply the following statement, which is similar in spirit
to the analysis in \S10 of~\cite\SHIMURAHECKE.

\proclaim{Conjecture 1} Let $A$ be an elliptic curve over~$\C$.
Then $A$
is modular
if and only if $A$ is isogenous
to all elliptic curves $\ssigma A$ obtained by conjugating
$A$ by algebraic automorphisms $\sigma$ of the field~$\C$.
\endproclaim

\noindent
Conjecture~1 is related to Serre's conjectures
via
a second conjecture, which we will state after
some motivating remarks.
Consider a normalized eigenform $f=\sum a_nq^n$
in the space of weight-two cusp forms on~$\Gamma_1(N)$.
Let $E=\Q(\ldots, a_n,\ldots)$ be the field generated by
the coefficients of~$f$, and let $d=[E\:\Q]$ be the 
degree of~$E$, i.e., the
dimension
of~$E$  as a $\Q$-vector space.
The abelian variety $A_f$ is an abelian variety over~$\Q$
of dimension~$d$ which
comes equipped with an action of~$E$.  To give sense to
the last statement,
one introduces the ring $\End(A_f)$ whose elements are 
maps $A_f\to A_f$ in the category of algebraic varieties over~$\Q$
which respect the group structure on~$A_f$.  (Such maps are
the endomorphisms of~$A_f$.)
The ring $\End(A_f)$ turns out to be
a free rank-$d$ module over~$\Z$; the
associated $\Q$-algebra 
$\End(A_f)\otimes\Q$
is isomorphic to~$E$.  Thus $A_f$ has many endomorphisms;
moreover, it is ``modular''
in the sense that there is a non-constant map
$X_1(N)\to A$ which is defined over~$\Q$.

\proclaim{Conjecture 2}
Let $A$ be an abelian variety over~$\Q$ for which
$\End(A)\otimes\Q$ is a number field of degree
equal to $\dim A$.
Then for some $N\ge1$,
there is
a non-constant map $X_1(N)\to A$ which is defined over~$\Q$.
\endproclaim

\noindent
In~\cite\KOREA, I prove that Serre's conjectures imply
Conjecture~2 and that Conjecture~2 implies Conjecture~1.

It is natural to regard Conjecture~1 and Conjecture~2 as
generalizations of the \Weil\ conjecture.
The first conjecture pertains to elliptic curves which are
not necessarily defined over~$\Q$, while the second pertains to
abelian varieties over~$\Q$ which are not necessarily
elliptic curves.  Neither of these conjectures
is proved in~\cite\WILES.
The \Weil\ conjecture can be generalized still further.
Indeed, a common generalization of
Conjectures
1 and~2 will be presented in a forthcoming work of E.~Pyle.
See also~\cite\FONTAINEMAZUR. 

\section Galois representations attached to elliptic curves\endhead
Let $A$ be an elliptic curve over~$\Q$.
A model for $A(K)$ when $K=\C$
is given by the Weierstra\ss\ theory
of complex analysis: the group $A(\C)$ is
the complex torus~$\C/L$,
where $L$ is the lattice of periods associated to the given
cubic equation.  (Explicitly, $L$ is
obtained by integrating the differential 
${\displaystyle dx\over\displaystyle y}$ 
on~$A$ over the
free abelian group $H_1(A(\C),\Z)$ of rank~two.)  Let $n$
be a positive integer, and let $A[n]$ be the group of
elements of~$A(\C)$ whose order divides~$n$.  This group
of {\it $n$-division points\/} on~$A$ may be modeled as 
${1\over n}L/L$; it is therefore a free module of rank two
over $\U n$, since $L$ is free of rank two over~$\Z$.

On reflection, one sees that $A[n]$ in fact lies in~$A(\Qbar)$,
where $\Qbar$ is the subfield of~$\C$ consisting of
all algebraic numbers.  Indeed, 
the group $A[n]$ is a finite subgroup of~$A(\C)$
which consists of those points satisfying
a certain set
of polynomial equations with rational coefficients; it follows
that the coordinates are algebraic numbers.
Moreover, let
$\GalQ$ be the group of automorphisms
of~$\Qbar$.  Then the same reasoning shows that
$A[n]$ is stable under the action of~$\GalQ$
on~$A(\Qbar)$ which results from the action of~$\GalQ$
on~$\Qbar$.  Thus
$A[n]$ comes equipped with a canonical action of the
Galois group $\GalQ$.

It is important to observe, for each $\sigma\in\GalQ$, that
the automorphism $P\mapsto\ssigma P$ is a {\it group\/} automorphism
of~$A[n]$;
in symbols,
$\ssigma(P+Q)=\ssigma P+\ssigma Q$ for
$P,Q\in A[n]$. 
This equation is a consequence of
the fact that the composition law which
expresses $P+Q$ in terms of $P$ and~$Q$ involves only
polynomials
with rational coefficients.
Our observation means that
the action of~$\GalQ$ on~$A[n]$ may be viewed
as a (continuous) homomorphism
$$ \rho_{A,n}\:\GalQ\to\Aut(A[n])$$
in which $\Aut(A[n])$ stands for the group
of automorphisms of~$A[n]$ as an {\it abelian group}.
Since $A[n]$ is isomorphic
to the group $(\U n)^2$, 
one has
$$\Aut(A[n])\approx \GL{\U n},$$
where the group on the right consists of 
two-by-two invertible matrices with coefficients
in~$\U n$.  While there is no canonical isomorphism
between these groups, each choice of basis
$A[n]\approx (\U n)^2$ determines such an isomorphism;
moreover,
the various isomorphisms obtained in this way
differ by inner automorphisms of~$\GL{\U n}$.  Therefore,
each element of $\Aut(A[n])$ has a well-defined trace
and determinant in~$\U n$.  

It is often fruitful to fix a choice of basis 
$A[n]\approx (\U n)^2$ and to view $\rho_{A,n}$
as taking values in the matrix group $\GL{\U n}$.
Once this choice is made, $\rho_{A,n}$ becomes 
matrix-valued
{\it representation\/} of the Galois group~$\GalQ$;
it is the representation of~$\GalQ$ defined by the
group of $n$-division points of~$A$.
The kernel of this representation corresponds, via Galois
theory, to a finite Galois extension $K_n$ of~$\Q$ in~$\Qbar$.
Concretely, this extension is obtained by
adjoining to~$\Q$ the coordinates of the
various points in~$A[n]$.  The Galois group $G_n\seteq\Gal(K_n/\Q)$
may thus be identified with
the image of $\rho_{A,n}$, which is a subgroup of
the target group $\GL{\U n}$.  The elliptic curve $A$
and the positive integer $n$ have given rise to a finite
Galois extension $K_n/\Q$ whose Galois group is a subgroup
of the group of two-by-two invertible matrices with
coefficients in~$\U n$.

It is natural to ask for a description of~$G_n$ as a
subgroup of~$\GL{\U n}$, cf.~\cite\SHCRELLE.
There is a 
(relatively rare) special case to
consider: that where $A$ has {\it complex multiplication\/} (over~$\C$).
When
$A$ is viewed as~$\C/L$,
the complex multiplication
case is that for which there is a complex
number $\alpha\not\in\Z$ such that $\alpha L\subseteq L$.
The group $G_n$ then has an abelian subgroup
of index~$\le2$, so it is much smaller than the ambient group
$\GL{\U n}$.  In the more common case where $A$ has no
complex multiplication, Serre showed 
in~\cite\SERREIM\ that the index of~$G_n$
in~$\GL{\U n}$ is bounded as a function of~$n$.
In particular,
$G_p = \GL{\U p}$ for all but finitely many primes~$p$.

As background, we point out that the \Weil\
conjecture was proved for complex multiplication elliptic
curves over~$\Q$ by Shimura in 
1971~\cite\SHIMURAPINK.  
The result of~\cite\SHIMURAPINK\ is suggestive and
may be regarded as evidence for the general case of the
\Weil\ conjecture.
As we recall below, however,
the elliptic curves which appear in
connection with Fermat's Last Theorem are semistable.
And it is a fact that
semistable elliptic curves
over~$\Q$ {\it never\/} have complex multiplication.
(One possible proof can be summarized as follows:
An elliptic curve with complex multiplication has an
integral $j$-invariant,
i.e., potentially good reduction.  Hence if it is
semistable, it has good reduction everywhere.
However, a theorem of~Tate
states that there is no elliptic curve over~$\Q$ with
everywhere good reduction, cf.~\cite\FONT.)
Accordingly, the theorem of~\cite\SHIMURAPINK\
cannot be used to prove Fermat's Last Theorem.

A key piece of information
about the extension $K_n/\Q$ (which depends on~$A$
as well as on~$n$) is that its discriminant
is divisible only by those prime numbers
which divide either~$n$ or the conductor of~$A$.
In other words, if $p\notdivide n$ is a prime number at
which $A$ has good reduction, then $K_n/\Q$ is unramified
at~$p$; one
says frequently that the representation
$\rho_{A,n}$ is unramified at~$p$.
Whenever this occurs, a familiar construction in
algebraic number theory produces
a Frobenius element
$\sigma_p$ in~$G_n$ which is well defined up to conjugation.

We shall now summarize this construction with $K_n$ replaced
by an arbitrary
finite Galois extension $K$ of~$\Q$.
Let 
${\Cal O}_K$
be the ring of algebraic integers in~$K$.
The Galois group $\Gal(K/\Q)$
leaves ${\Cal O}_K$ invariant, so that one obtains an
induced action of~$\Gal(K/\Q)$ on the ideals of~${\Cal O}_K$.
The
set of
prime ideals $\goth p$ of~${\Cal O}_K$ which contain the prime
number~$p$ is permuted under this action.
For each~$\goth p$, the
subgroup $D_{\goth p}$
of~$\Gal(K/\Q)$ consisting of those elements in~$\Gal(K/\Q)$ which
fix $\goth p$ is called the {\it decomposition group\/} of~$\goth p$
in~$\Gal(K/\Q)$.  Meanwhile, the finite field $\F_{\goth p}\seteq
{\Cal O}_K/\goth p$ is a finite extension of the prime field $\Fp$.
The extension $\F_{\goth p}/\Fp$ is necessarily Galois;
its Galois group is the cyclic group generated
by the Frobenius automorphism $$\phi_p\: x\mapsto x^p$$
of~$\F_{\goth p}$.  There is a natural map 
$D_{\goth p}\to\Gal(\F_{\goth p}/\Fp)$, gotten by associating a
given $\delta\in D_{\goth p}$
to the automorphism of ${\Cal O}_K/\goth p$
induced by~$\delta$.  This map is surjective;
its injectivity is equivalent to the assertion that $p$ is
unramified in the extension $K/\Q$.  
Therefore, whenever this assertion is true,
there is a unique
$\sigma_{\goth p}\in D_{\goth p}$ 
whose image in~$\Gal(\F_{\goth p}/\Fp)$ is~$\phi_p$.
The automorphism $\sigma_{\goth p}$ is then a well defined
element of~$\Gal(K/\Q)$, known as the Frobenius automorphism
for~$\goth p$. It is easy to show that the various
$\goth p$ are all conjugate under~$\Gal(K/\Q)$ and that the
Frobenius automorphism for the conjugate of~${\goth p}$ by~$g$
is~$g\sigma_{\goth p}g^{-1}$.  In particular, the various
$\sigma_{\goth p}$ are all conjugate; this justifies the
practice of writing $\sigma_p$ for any one of them and
stating that $\sigma_p$ is well defined up to conjugation.

For later use, we prolong this digression and introduce
the concept of Frobenius elements in~$\GalQ$.
Let $p$ again be a prime number, and let $\goth p$ now be
a prime of~$\Qbar$ lying over~$p$.  (One can think of $\goth p$
as a coherent set of choices of primes lying over~$p$ in the
rings of integers of all finite extensions of~$\Q$ in~$\Qbar$.)
To~$\goth p$ we associate: (1) its residue field $\F_{\goth p}$,
which is an algebraic closure of the finite field $\F_p$,
and (2) a decomposition subgroup $D_{\goth p}$ of~$\GalQ$.
There is again a surjective map
$D_{\goth p} \to \Gal(\F_{\goth p}/\F_p)$.
The Frobenius
automorphism $\phi_p\: x\mapsto x^p$
generates the
target group
in the topological sense:
the subgroup of~$\Gal(\F_{\goth p}/\F_p)$ consisting of
powers of~$\phi_p$
is dense in~$\Gal(\F_{\goth p}/\F_p)$.  
We shall use the
symbol $\Frob_p$ to denote any preimage of 
$\phi_p$ in~$D_{\goth p}$ and refer to~$\Frob_p$ as a Frobenius
element for~$p$ in~$\GalQ$.

One thinks of $\Frob_p$ as a specific element of~$\GalQ$,
albeit one which is doubly ill-defined.  The ambiguities
in~$\Frob_p$ result from the circumstance that $\goth p$ needs to
be chosen and from the fact that 
$D_{\goth p}\to\Gal(\F_{\goth p}/\F_p)$ has a large kernel,
the inertia subgroup $I_{\goth p}$
of~$D_{\goth p}$.
The usefulness of $\Frob_p$ stems from the fact that
the various $\goth p$ are conjugate by elements of~$\GalQ$,
so that all the subgroups $D_{\goth p}$ of~$\GalQ$ are
conjugate, and similarly all $I_{\goth p}$ are conjugate.
Thus
if $\rho$ is a homomorphism mapping $\GalQ$ to some
other group, the kernel of
$\rho$ contains one $I_{\goth p}$ if and only if
it contains all $I_{\goth p}$.  In this case, one says
that $\rho$ is unramified at~$p$; the image of~$\Frob_p$
is then an element of the target group of~$\rho$ which
is well defined up to conjugation.
Note that one may write $\rho_{A,n}(\Frob_p)=\sigma_p$ for
all primes $p$ at which $\rho_{A,n}$ is unramified.  As has
been stated, these include all primes which divide neither~$n$
nor the conductor of~$A$.

Let us return now to the matrix group~$G_n$.
We pointed out above that each element of~$G_n$ has a
trace and determinant in~$\U n$ which are independent of
any choice of basis.
The Frobenius element
$\sigma_p$ is an element of~$G_n$
which is well defined only up to conjugation.
Nevertheless, the trace and determinant of~$\sigma_p$
are well-defined, since conjugate matrices have the same
traces and determinant.
The number
$\det \sigma_p$ is the residue class of $p$ mod~$n$.  On the
other hand, one has the striking congruence
$$ \tr(\sigma_p) \equiv b_p \hbox{ mod }n,$$
where $b_p$ is the number $p{+}1-\#(A(\Z/p\Z))$ introduced
above.  This means that the representation $\rho_{A,n}$
encapsulates information about the numbers $b_p$ (for 
prime numbers $p$ which are primes of good reduction and which
are prime to~$n$); more precisely, it determines the
numbers $b_p$ mod~$n$.

\section Galois representations attached to modular forms\endhead
Suppose that $f\in S(N)$ is a normalized eigenform.  
The coefficients $a_n$ of~$f$ are always algebraic integers,
but not necessarily ordinary integers.
If it happens that the $a_n$ all lie in~$\Z$, then
the abelian variety $A=A_f$ is an elliptic curve.
By considering the
family $\rho_{A,n}$, one obtains a series of
representations of the Galois group~$\GalQ$.  These representations
are related to~$f$ by the congruence 
$ \tr(\sigma_p) \equiv a_p \hbox{ mod }n$, valid for the $n$th
representation and all primes $p\notdivide nN$.  We are especially
interested in the case where $n$ is a prime number~$\ell$;
the ring $\U n$ is then the {\it finite field}~$\F_\ell$.

The representations
$\rho$ are associated to~$A$, which in turn arises from~$f$.
Hence it is tempting to write $\rho_{f,\ell}$
for the representations $\rho_{A,\ell}$.  The obstacle to
doing this arises from the circumstance
that $A$ is determined
up to isogeny, but not always
up to isomorphism.  If one replaces $A$
by an isogenous elliptic curve, the representations $\rho_{A,\ell}$
may change!  To circumvent this difficulty, we introduce the
``semisimplifications'' of the~$\rho_{A,\ell}$.

These representations are defined as follows.
If $\rho$ is a two-dimensional
representation of a group over a field, $\rho$ is either
irreducible, or else ``upper-triangular,'' i.e., an extension
of a one-dimensional representation
$\alpha$ by another, $\beta$.  In the case
where $\rho$ is irreducible, we declare its semisimplification
to be $\rho$ itself.  In the reducible case, the semisimplification
of~$\rho$ is the {\it direct sum\/} of the two one-dimensional
representations $\alpha$ and~$\beta$.
Clearly, the trace and determinant are the same
for~$\rho$ and for its semisimplification.

When $A$ is fixed, the results of~\cite\SERREIM\
show that $\rho_{A,\ell}$ can be reducible only
for a finite number of~$\ell$.  In fact, a theorem of
Mazur~\cite\PRIMEDEGREE\ shows that $\rho_{A,\ell}$
is irreducible for all $\ell$ not in the set
$\{\,2, 3, 5, 7, 13, 11, 17, 19, 37, 43, 67, 163\,\}$.  
Hence the replacement
of~$\rho_{A,\ell}$ by its semisimplification can be
thought of as ``fine tuning'' which affects only a
small number of the representations.
One shows
easily for all~$\ell$
that the semisimplification of~$\rho_{A,\ell}$ 
depends only on~$f$ and on~$\ell$ (but not on the choice of~$A$).
Introducing
$$ \rho_{f,\ell} \seteq \hbox{semisimplification of }\rho_{A,\ell},$$
one obtains a sequence of semisimple representations of~$\GalQ$
which are well defined up to isomorphism.  The characteristic
property of~$\rho_{f,\ell}$ may be summarized in terms of
Frobenius elements $\Frob_p$ in the Galois group $\GalQ$,
elements which were introduced above.
Namely, if
$p$ is a prime number not dividing $\ell N$, then
$\rho_{f,\ell}(\Frob_p)$ has trace $a_p$ mod~$\ell$ and determinant
$p$ mod~$\ell$.

It is natural to generalize this process by considering the
situation where 
$f=\sum a_n q^n$ is a
normalized eigenform whose coefficients $a_n$
are 
algebraic integers, but
not necessarily
rational integers.
As we indicated above, the field $E=E_f$
generated by the~$a_n$ is a number field, i.e., a
finite extension of~$\Q$.  Moreover, the coefficients $a_n$
of~$f$ lie in the integer ring ${\Cal O}_f$ of~$E$.
It is perhaps worth noting that the ring $\Z[\ldots,a_n,\ldots]$
generated by
the $a_n$ inside~$E$, while a subring of finite index
in~${\Cal O}_f$, is not necessarily equal to~${\Cal O}_f$.

Using the abelian variety $A_f$, one constructs
representations indexed not by the prime {\it numbers}, but
rather by the non-zero prime {\it ideals\/} of~${\Cal O}_f$.
(For details, see~\cite{\SHIMURABOOK, Ch.~7}.)
If $\lambda$ is such a prime, its residue field $\F_\lambda$ is
a finite field, say of characteristic~$\ell$.  The prime
field $\F_\ell=\U\ell$ is then canonically embedded in~$\F_\lambda$.
For each $\lambda$, one finds a semisimple representation
$\rho_{f,\lambda}\:\GalQ\to\GL{\F_\lambda}$
which is characterized up to isomorphism
by the following property: if $p$ is
a prime number not dividing $\ell N$, then 
$\rho_{f,\lambda}(\Frob_p)$ has trace $a_p$ mod~$\lambda$ and
determinant $p$ mod~$\lambda$.

The assertion concerning the determinant 
of the matrices~$\rho_{f,\lambda}(\Frob_p)$ 
may be rephrased as the statement that
the
determinant of
the {\it representation\/} $\rho_{f,\lambda}$
is the mod~$\ell$ cyclotomic character~$\chi_\ell$.
This character is defined by considering
the group $\mu_\ell$ of $\ell$th roots of unity in~$\Qbar$;
the action of the
Galois group $\GalQ$ on the cyclic group~$\mu_\ell$
gives rise to a continuous homomorphism
$$\chi_\ell\: \GalQ\to\Aut(\mu_\ell).$$
Since $\mu_\ell$ is a cyclic group of order~$\ell$,
its group of automorphisms is canonically the group 
$\UU\ell=\F_\ell^\ast$.  We emerge with a 
map $\GalQ\to\F_\ell^\ast$,
which is the character in question.
The equality $$\det\rho_{f,\lambda}
=\chi_\ell$$
interpreted by viewing
both homomorphisms as taking values in~$\F^\ast_\lambda$.

Suppose now that $c\in\GalQ$ is the
automorphism ``complex conjugation.''  Then the determinant
of~$\rho_{f,\lambda}(c)$ is~$\chi_\ell(c)$.  Now $c$
operates on roots of unity by the map $\zeta\mapsto\zeta^{-1}$,
since roots of unity have absolute value~1.  Accordingly,
$$\det(\rho_{f,\lambda}(c))=-1;$$ one says that $\rho_{f,\lambda}$
is {\it odd}. 

This parity statement generalizes to modular forms ``with
Nebentypus'' whose weights are not necessarily two.
Here is a quick synopsis of the
situation; some relevant references are provided
in~\cite\RIBETSURVEY.  For integers $k\ge1$,
$N\ge1$ and characters $\epsilon\:\UU N\to\C^\ast$, one
considers the space $S_k(N,\epsilon)$ of weight-$k$ cusp
forms with character $\epsilon$ on~$\go N$; we have
$S(N)=S_2(N,1)$.  
This space is automatically zero unless $\epsilon(-1)=(-1)^k$,
so we will assume that this parity condition is satisfied.
The space $S_k(N,\epsilon)$ admits an
operation of Hecke operators $T_n$, and we again have the
concept of a normalized eigenform in~$S_k(N,\epsilon)$.
If $f=\sum a_nq^n$ is such a form, the numbers $a_n$ ($n\ge1$)
and the values of~$\epsilon$ all lie in a single integer 
ring~${\Cal O}_f$. For each non-zero prime ideal $\lambda$ of~${\Cal O}_f$,
one constructs a semisimple representation
$$\rho_{f,\lambda}\:\GalQ\to\GL{\F_\lambda}.$$
Let $\ell$ again denote the characteristic of~$\F_\lambda$.
Then for all $p\notdivide \ell N$, the trace 
of~$\rho_{f,\lambda}(\Frob_p)$ is again $a_p$ mod~$\lambda$.
The determinant of this matrix is 
$p^{k-1}\epsilon(p)$ mod~$\lambda$.

Once the proper definition is made,
the determinant of the map $\rho_{f,\lambda}$
becomes the product $\chi_\ell^{k-1}\epsilon$.
In writing $\det\rho_{f,\lambda}=\chi_\ell^{k-1}\epsilon$,
we use $\chi_\ell$ to denote the mod~$\ell$ cyclotomic
character and employ a standard construction to regard
$\epsilon$
as a map $\GalQ\to\F_\lambda^\ast$.
The construction in question begins with the
map $\GalQ\to\UU N$ giving the action of~$\GalQ$
on the $N$th roots of unity.  Composing this map
with the character $\epsilon$, we obtain
a homomorphism $\GalQ\to{\Cal O}_f^\ast$.  On reducing
this homomorphism mod~$\lambda$, we obtain the desired
variant of~$\epsilon$.

Evaluating the formula 
$\det\rho_{f,\lambda}=\chi_\ell^{k-1}\epsilon$
on $c=\hbox{``complex
conjugation''}$,
one
finds
$$\det(\rho_{f,\lambda}(c)) = (-1)^{k-1}\epsilon(c) = -1.$$
In these equalities, we exploit the fact that
$\epsilon(c)$ is another name for~$\epsilon(-1)$ and
remember the parity condition
$\epsilon(-1)=(-1)^k$.  The upshot of this
is that the representations
$\rho_{f,\lambda}$ are always odd, even in the generalized
set-up.

In fact, it is possible in this situation to find a
normalized eigenform
$f'$ of weight~two with some character~$\epsilon'$, along with
a maximal ideal $\lambda'$ of the integer ring for~$f'$,
so that the representations $\rho_{f,\lambda}$
and~$\rho_{f',\lambda'}$ are isomorphic,
cf.~\cite{\MOTIVES, Th.~2.2 and Cor.~3.2}.  (To compare
these representations, it is necessary to embed the residue fields
of~$\lambda$ and of~$\lambda'$ in a 
suitably chosen
common finite field
of characteristic~$\ell$.)  Hence the ``generalized set-up''
can be reduced to the case of weight~two, provided that one
considers eigenforms
for which the associated characters
may be non-trivial.  The process of reduction to weight two
is a very powerful one in the theory, since the representations
arising from forms of this weight are constructed directly
from points of finite order on abelian varieties.
To the best of my knowledge, the idea 
of reducing systematically to weight two
originated with an
unpublished 1968 manuscript of Shimura~\cite\SHIMURAOLD.

Before leaving this topic, we should mention that the cases
$k=1$ and $k>1$ are quite distinct in flavor.  In the former
case, the representations $\rho_{f,\lambda}$
may be viewed as the set of
reductions of a single continuous
representation $\rho_f$
with finite image.
For details, see~\cite\DELIGNESERRE.

\section Serre's conjectures\endhead
We shall give a brief summary of conjectures made by
J.-P.~Serre in~\cite\DUKE.
A recent article
by~H.~Darmon~\cite\NEWDARMON\ 
discusses the conjectures more extensively and emphasizes
applications and numerical examples.

Let $\ell$ be
a prime number
and let $\F$ be an algebraic closure of the prime field
$\F_\ell$.  Suppose that $\rho\:\GalQ\to\GL\F$ is an odd
continuous representation.
It seems natural to ask whether or not $\rho$ is 
the mod~$\lambda$ representation attached to a suitable
normalized
eigenform.  Since the representations $\rho_{f,\lambda}$
are semisimple by definition, it is necessary to
limit our discussion to the case where $\rho$ is semisimple.
In fact, the case where $\rho$ is semisimple and reducible
is sometimes awkward, so we will assume from now on that
$\rho$ is irreducible.  However, the excluded
case where $\rho$ is reducible is quite interesting; see \cite\BUZ\
for some observations in this case.

Let us say then that
$\rho$ is {\it modular\/} if
one can find: (i) a normalized eigenform $f$ in
some space $S_k(N,\epsilon)$; (ii) a prime $\lambda$
dividing $\ell$ in the ring of integers ${\Cal O}_f$ associated
to~$f$; and (iii) an embedding $\F_\lambda\hookrightarrow\F$
such that $\rho$ is isomorphic to the representation
obtained by composing $\rho_{f,\lambda}$ with the
inclusion $$\GL{\F_\lambda}\hookrightarrow\GL\F$$ associated
with~(iii).  A weak form of the conjectures made by~Serre
in~\cite\DUKE\ is the following statement: {\it Every irreducible
continuous odd representation
$\rho\:\GalQ\to\GL\F$ is modular}.

This statement
was first formulated by~Serre in the 1970s
for modular forms of level~1 (i.e., on~$\SL\Z$).  For such forms,
the representations $\rho_{f,\lambda}$ are ramified only
at~$\ell$; Serre asked whether an odd irreducible $\rho$
which is unramified outside~$\ell$ is necessarily
associated to a normalized eigenform on~$\SL\Z$.
Tate confirmed this for $\ell=2$ \cite\TATELETTER\ by showing
that there are no such representations.
(As Tate remarks at the end of his article, Serre treated
the case $\ell=3$ in a similar manner by exploiting the
discriminant bounds of Odlyzko and Poitou.)
The case where $\rho$ may be ramified at primes other than~$\ell$
was taken up by Serre in
the mid-1980s, when
computer calculations
by J.-F.~Mestre convinced Serre
that the conjectured statement was plausible.

The qualitative conjecture
to this effect is supplemented in~\cite\DUKE\ by
an intricate recipe which pinpoints the space $S_k(N,\epsilon)$
where one should find an eigenform~$f$ giving rise to
a specific representation~$\rho$.  (As Serre later
observed, the recipe for~$\epsilon$
needs to be modified in certain cases when $\ell=2$ or~3.)
The conjunction of the qualitative statement 
that $\rho$ is modular
and the precise
recipe 
fingering the space $S_k(N,\epsilon)$ 
is sometimes called the Strong Serre Conjecture.
One possible justification for this name is the fact that
the ``strong'' conjecture
immediately implies
the \Weil\
conjecture, Fermat's Last Theorem, and a host of other
assertions! 
It has gradually emerged that the qualitative statement and
its strong cousin are in fact {\it equivalent}, at least
when $\ell>2$; see \cite\MOTIVES\ and \cite\DIAMONDREF\
for a proof of the equivalence.
Hence it is now possible
to use the singular term ``Serre's conjecture'' to
refer to what was initially a package of interrelated conjectures.

Perhaps I should close this section by expressing
the sentiment that
the conjecture of Serre, while visibly important,
currently seems intractable.  Given an irreducible
representation $\rho$ with odd determinant, one is at a loss
for a strategy which will lead to a proof that $\rho$
is modular.  In my Progress in Mathematics
lecture, I engaged in
a certain amount of philosophical
speculation, stressing
the parallel between the \Weil\
conjecture and Serre's conjecture.  Each conjecture
states that all
objects of a certain type are modular; in both cases one
has a unidirectional ``arrow''\emrule
a means of constructing objects from modular
forms.  Thus one's impulse is to try to
attack these conjectures by an appropriate form of
``counting.''  While Serre's conjecture is broader than the
\Weil\
conjecture (the former implies the latter), one might
suspect
that Galois representations might be easier to count than
elliptic curves.  To the extent that this is so,
one could imagine attacking the geometric conjecture
about elliptic curves via
the Galois-theoretic
conjecture of~Serre.  

This philosophy
is perhaps not far removed in spirit from the strategy used
by Andrew Wiles in
approaching the
\Weil\
conjecture.
Certainly, however, the analogy should not be taken too
seriously; in fact, Wiles introduced his approach in~1993
with the statement that
it was ``orthogonal'' to Serre's conjecture.

\section Frey's construction\endhead
Like most recent 
work on Fermat's Last Theorem, the connection
between Serre's conjecture and~FLT begins with
constructions linking solutions to Fermat's equation
with elliptic curves.  Although Y.~Hellegouarch and
others had noted such constructions,
a decisive step
was taken by G.~Frey in 
an unpublished 1985 manuscript entitled ``Modular
elliptic curves and Fermat's conjecture.''

Frey's idea goes as follows.
Suppose that there is a non-trivial
solution to Fermat's equation
$X^\ell+Y^\ell=Z^\ell$.  We can assume that the exponent
is
a prime number different from 2 and~3 and that
the solution is given by a triple of 
relatively prime integers $a$, $b$, and~$c$.
The equation
$y^2=x(x-a^\ell)(x+b^\ell)$ then
defines an elliptic curve $A$
with unexpected properties.

These properties are catalogued in~\S4.1 of Serre's
article~\cite\DUKE:
After performing some elementary manipulations, we arrive
at a triple $(a,b,c)$ in which $b$ is even and $c$ is
congruent to~1 mod~4.  Frey's construction yields for this
triple an elliptic curve~$A$ whose conductor
is the product of the prime
numbers which divide $abc$
(each occurring to the first power).  In particular,
$A$ is semistable.
On the other hand, the minimal discriminant
$\Delta$ of~$A$
is the quotient of~$(abc)^{2\ell}$ by the
factor $2^8$.  From Frey's point of view,
the main ``unexpected'' property of~$A$ is that
$\Delta$ is the product of a 
power of~2 and a
perfect $\ell$th power, where $\ell$ is a prime~$\ge5$.
Frey translated this property into a statement about the N\'eron
model for~$A$: if $p$ is an odd prime at which $A$ has bad
reduction, the number of components in the mod~$p$ reduction
of the N\'eron model is divisible by~$\ell$.  Frey's idea
was to compare this number to
the corresponding
number for the Jacobian of the modular curve $\xo N$, where $N$
is the conductor of~$A$.  Frey predicted that
a discrepancy between the two
numbers would preclude $A$ from being modular.
In other
words,
Frey concluded
heuristically that the existence
of~$A$ was incompatible
with the \Weil\ conjecture, which asserts that
all elliptic curves over~$\Q$ are modular.

Frey's construction spawned several lines of inquiry,
in which mathematicians sought either to prove Fermat's
Last Theorem outright, or else
to link it to established or emerging conjectures
such as the $abc$ conjecture and
Szpiro's conjecture.
These
latter conjectures
are
treated  by
such articles as
\cite\FREYONE,
\cite\FREYTWO,
\cite\HINDRY,
\cite\LANG,
\cite\MAZURNUMBER\
and~\cite\OESTERLE. From
the point of view of
Szpiro's conjecture and
the $abc$ conjecture, the surprising feature
of Frey's curve is 
the size of its discriminant, 
rather than any special properties of the
discriminant's factorization.  For Frey's curve $A$,
the absolute value of the discriminant can be bounded
from below by a
high power of the conductor of~$A$.

To illustrate the force of Frey's construction, we
now sketch the deduction of Fermat's Last Theorem
from Serre's conjectures~\cite\DUKE.
(A word of caution: these
conjectures are still conjectures!)
Suppose
that $a$, $b$ and~$c$ 
are relatively prime integers which
satisfy Fermat's equation with exponent~$\ell$.  After
performing the manipulations mentioned above, we may
write down a Frey curve $A$ whose discriminant is the product of a 
power of~2 and a
perfect $\ell$th power.  If $\ell$ is greater
than~3, a theorem of Mazur~\cite\PRIMEDEGREE\ asserts that the
representation $\rho=\rho_{A,\ell}$ defined by~$A[\ell]$
is an irreducible representation of~$\GalQ$.
Because of the hypothesis on the discriminant of~$A$, $\rho$
is unramified outside $2$ and~$\ell$; moreover it is
``finite'' at~$\ell$ in a sense which is explained in~\cite\DUKE.
The recipes of Serre's article
require that
$\rho$ arises from a
normalized eigenform in the space~$S(2)$.
However, as was mentioned above,
this space has dimension~0.

\section Conjecture ``epsilon''\endhead
A first step toward justifying
Frey's heuristic conclusion
was taken in August, 1985 by
Serre in a letter to J.-F.~Mestre~\cite\ARCATA.
In this letter,
the writer formulated two related conjectures about
modular forms, which he called
$C_1$ and~$C_2$.  
(These conjectures predate the conjectures of~\cite\DUKE;
they are now special cases of the latter conjectures.)
Serre pointed out
that Fermat's Last Theorem is a consequence
of the \Weil\ conjecture {\it together with\/}
the two new conjectures.
As it was thought initially that $C_1$ and~$C_2$ would be
easy to establish, the two statements
immediately acquired the collective
nickname ``Conjecture~$\epsilon$.''
One thus had
$$ \hbox{\Weil} +\epsilon
    \,\Longrightarrow\,\hbox{Fermat's Last Theorem}.$$

Serre's ``$\epsilon$'' conjecture is in fact a subtle statement about
the
mod~$\ell$ Galois representations 
arising from eigenforms
in the various spaces~$S(N)$. Specifically, suppose that $f\in S(N)$
is a normalized eigenform, and let $\lambda$ be a prime ideal in
the ring of integers of the field generated by the coefficients of~$f$.
Then $\rho_{f,\lambda}$ is a semisimple representation of~$\GalQ$
with values in a finite field, whose characteristic we will call~$\ell$.
This representation is unramified at all primes $p\ne\ell$ which
do not divide~$N$; in other words, $\rho_{f,\lambda}$ has the
right to be ramified at~$\ell$ and at those primes
which divide~$N$.  Serre's conjecture concerns the case where
$\rho_{f,\lambda}$ is unramified at such a prime: it predicts
that this behavior can be attributed to the
existence of
a normalized eigenform $f'$ of level lower than~$N$
which gives rise to~$\rho_{f,\lambda}$.

Specifically, suppose that $\rho_{f,\lambda}$ is irreducible,
that $\ell$ is different from~2, and that $\rho_{f,\lambda}$
is unramified at a prime number $p\ne\ell$ which divides~$N$ but whose
square does not divide~$N$.  Serre's conjecture predicts that
there is an eigenform $f'\in S(N/p)$, together with a prime ideal
$\lambda'$ in the integer ring of the field of coefficients of~$f$,
such that $\rho_{f,\lambda}$ and $\rho_{f,\lambda'}$ are isomorphic.
A variant of this conjecture concerns the case $p=\ell$.
(Note that, as
in a situation discussed earlier,
one must embed the residue fields
of~$\lambda$ and of~$\lambda'$ in a 
suitably chosen
common finite field
of characteristic~$\ell$
before
comparing $\rho_{f,\lambda}$ and~$\rho_{f,\lambda'}$.)

I proved Serre's level-lowering conjecture
in a 1990 article~\cite\FERMAT, thereby
establishing
the
implication
  $$\hbox{Conjecture of \Weil\ }
       \Longrightarrow\hbox{ Fermat's Last Theorem}$$
which was the goal of Frey's construction.
(See also \cite\PRASAD\ and~\cite\TOULOUSE\
for expository accounts of this work.)
Since the Frey curves associated with Fermat solutions
are semistable elliptic curves, I proved that 
the semistable case of the \Weil\ conjecture implies
Fermat's
Last Theorem.

To illustrate the logic used in establishing this implication,
we consider a semi\-stable elliptic curve $A$ over~$\Q$,
together with a prime $\ell\ge3$ for which the representation
$\rho_{A,\ell}$ is irreducible.  Suppose that $\rho_{A,\ell}$
is unramified at all primes $p\ne2,\ell$ and moreover that
$\rho_{A,\ell}$ is ``finite'' at~$\ell$.  Then the results 
of~\cite\FERMAT\ assert that $A$ cannot be modular.  
To apply these results to Fermat's Last Theorem, we
suppose that $A$ is the
Frey curve associated to a hypothetical solution
to the degree-$\ell$ Fermat equation with $\ell>3$.
Then,
as was noted above,
$A$ is semistable,
and $\rho_{A,\ell}$ 
is an irreducible representation with the indicated
ramification and ``finiteness'' properties.
Accordingly,
$A$
cannot be modular.
Consequently, if
all 
elliptic
curves over~$\Q$ are modular, then there can be no solution to
Fermat's equation.

Further light can be shed on~\cite\FERMAT\ if we focus on
the simplest situation in which its results apply.
Suppose that $A$ is an elliptic curve over~$\Q$ whose conductor
is a prime number~$p$.  Let $\ell$ be a prime number different
from 2 and~$p$ for which $\rho_{A,\ell}$ is irreducible.
The representation $\rho_{A,\ell}$ is unramified at all primes
other than $p$ and~$\ell$, and the contribution of~\cite\FERMAT\
is to show that $\rho_{A,\ell}$ is indeed ramified at~$p$ if
$A$ is modular.  To prove this, the one supposes
that $A$ is modular, so that $\rho_{A,\ell}$
is connected up with the
space of weight-two cusp forms on~$\go p$.  It is possible to
show
that $\rho_{A,\ell}$ is similarly connected with other
discrete subgroups of~$\SL{\bold{R}}$, coming from
indefinite quaternion division algebras over~$\Q$.
More precisely, there are prime numbers $q\ne p$ such that
$\rho_{A,\ell}$ arises from the quaternion algebra over~$\Q$
of discriminant $pq$.  The desired result
about $\rho_{A,\ell}$ follows from a detailed comparison
of the mod~$p$ and mod~$q$
reductions
of the three modular curves $\xo p$, $\xo q$ and $\xo{pq}$
with the mod~$p$ and mod~$q$ reductions of
the Shimura curve associated with the quaternion algebra
of discriminant~$pq$.
The latter curve is an analogue of~$\xo{pq}$ in which
the group $\go{pq}$ is replaced by the group of norm-1
elements in a maximal order of the quaternion algebra of
discriminant~$pq$.

\section Wiles's strategy\endhead
Suppose that $A$ is an elliptic curve over~$\Q$.
To verify the \Weil\ conjecture for~$A$
is to link $A$ to modular forms.  In an approach
inspired by Serre's conjecture, one
might
begin by
considering
the representations $\rho_{A,\ell}$ obtained
from the action of~$\GalQ$ on~$A[\ell]$, when $\ell$ is a prime
number. If one could show that an infinite number of these
representations are modular (in the broadest possible
sense), one would go on to prove that $A$ is modular.
Alas, as was indicated above,
it is not clear how to translate this approach into a proof.

We mentioned previously that Wiles's approach to the
\Weil\ conjecture is ``orthogonal'' to one
based on consideration of the varying~$\rho_{A,\ell}$.
Here is the nub of the idea:
One first fixes a prime~$\ell$ and considers the
family of groups $A[\ell^\nu]$ for $\nu=1,2,\ldots$. The
resulting sequence of representations 
$$ \rho_{A,\ell^\nu}\:\GalQ\to\GL{\U{\ell^\nu}}$$
may be packaged as a single representation
$$ \rho_{A,\ell^\infty}\:\GalQ\to\GL{\Z_\ell},$$
where $\Z_\ell$ is the ring of $\ell$-adic integers,
i.e., the projective limit of the rings $\U{\ell^\nu}$.
To
prove that $A$ is a modular elliptic curve, it suffices
to show that $\rho_{A,\ell^\infty}$ is modular in an
appropriate sense.  Indeed, the trace of 
$\rho_{A,\ell^\infty}(\Frob_p)$ coincides with the 
rational integer $b_p$
for all $p\notdivide \ell N$, where $N$ is the conductor of~$A$.
As soon as one finds an eigenform $f$ in~$S(N)$ whose
eigenvalues are related to the
traces of~$\rho_{A,\ell^\infty}(\Frob_p)$, one has essentially
proved that $A$ is modular.

Needless to say, 
if $\rho_{A,\ell^\infty}$ is modular, then so, in particular,
is~$\rho_{A,\ell}$.
Relating $\rho_{A,\ell^\infty}$ to
modular forms is thus no easier than
the formidable task of
proving that
$\rho_{A,\ell}$ is modular!
On the other hand, to prove that $A$ is modular by the
$\ell$-adic method, we need only work with a {\it single\/}
prime~$\ell$.  The approach of~\cite\WILES\
capitalizes
on the fact
that the 
finite groups $\GL{\F_2}$ and~$\GL{\F_3}$ are solvable.
This circumstance enables one to apply deep results
of Langlands~\cite\LANGLANDS\ and Tunnell~\cite\TUNNELL\
to prove that
$\rho_{A,\ell}$ is~modular
for $\ell\le3$, cf.~\cite{\RUBINSILVER, \S2.3}.
(The relevant results of Langlands are those
concerning
the theory of base change
\`a la Saito-Shintani. 
For expositions of these results, see
\cite\GELBART\ and \cite\LABESSE.)

Wiles's basic idea is to
prove that that if $\ell$ is a prime for which
$\rho_{A,\ell}$ is modular, then
$\rho_{A,\ell^\infty}$
is automatically modular (and hence $A$
is a modular elliptic curve).
In thinking about the jump from~$\rho_{A,\ell}$
to~$\rho_{A,\ell^\infty}$, ignores $A$ as
much as possible\emrule the aim is to prove results
about $\ell$-adic representations which can be applied
to~$\rho_{A,\ell^\infty}$.

\section The language of deformations\endhead
We now introduce the machinery which underlies
the jump from the modularity of~$\rho_{A,\ell}$ to the
modularity of~$\rho_{A,\ell^\infty}$.
We suppose for simplicity that
$A$ is a semistable elliptic curve over~$\Q$,
and we let
$\ell\ge3$ be a prime number
for which $\rho_{A,\ell}$
is both modular and~irreducible.  
Choosing a basis of~$A[\ell]$, we regard $\rho_{A,\ell}$
as taking values in the matrix group~$\GL{\F_\ell}$.
As was suggested above,
we seek to establish
the modularity of~$\rho_{A,\ell^\infty}$
by a method which treats simultaneously
{\it all\/}
lifts of $\rho_{A,\ell}$
which are plausibly modular.

In this context, lifts are continuous
homomorphisms $$\tilde\rho\:\GalQ\to\GL B,$$ where $B$
is a complete local Noetherian $\Z_\ell$-algebra with residue
field~$\F_\ell$.  They are constrained to lift~$\rho_{A,\ell}$
in the obvious sense.  Namely,
we require that
$\rho_{A,\ell}$ coincide with the
composite of~$\tilde\rho$ and the homomorphism
$\GL B\to\GL{\F_\ell}$ induced by
the residue map $B\to\F_\ell$.
(his depiction
ignores a technical wrinkle: it might to
necessary later on
to replace
$\Z_\ell$ by the integer ring of a finite extension of
the $\ell$-adic field~$\Q_\ell$.)
The lifts which are ``plausibly modular'' are
those which obey a set of local properties.
The word ``local'' is meant to allude to the restrictions of
$\tilde\rho$ to the subgroups 
$\Gal(\Qpbar/\Qp)$ of~$\GalQ$
obtained as decomposition groups for prime numbers~$p$
(or, more precisely, for primes of~$\Qbar$).  Wiles
imposes conditions on these restrictions which lift
conditions already satisfied by~$\rho_{A,\ell}$.  
These conditions are summarized in~\S3.4 of~\cite\RUBINSILVER;
we shall evoke them below as well.
There is
flexibility and tension
implicit in the choice
of these conditions.  They should be broad enough to be
satisfied by~$\rho_{A,\ell^\infty}$ and tight enough to
be satisfied only by lifts that can be related to modular forms.
Roughly speaking,
in order to prove the modularity of all lifts satisfying
a fixed set of conditions,
you need to
specify in advance a space
of modular forms $S$ so that the 
normalized
eigenforms in~$S$ satisfy the
conditions and such that, conversely, all lifts satisfying
the conditions are plausibly related to forms in~$S$.
It is intuitively clear that
this program will be simplest to
carry out when the conditions are the most stringent
and progressively harder to carry out as the conditions
are relaxed.

A theme which emerges rapidly is that there are at least two 
sets of conditions of special interest.
Firstly,
one is especially at ease when
dealing with the most stringent
possible set of conditions which are satisfied by~$\rho_{A,\ell}$;
this leads to what Wiles calls the ``minimal'' problem.
Secondly,
one needs at some point to consider
some
set of conditions which allows
treatment of the lift $\rho_{A,\ell^\infty}$\emrule
this lift is, after all, our main target.
It would be natural to consider the most stringent such set.
The two sets of conditions may coincide, but there is 
no guarantee that they do;
in general, the second set of conditions is more generous
than the first.

In~\cite\WILES, Wiles
provides a beautiful ``induction'' argument which enables him
to pass from the minimal set of conditions to a
non-minimal set.
Heuristically, this argument requires
keeping tabs on the set of those 
normalized
eigenforms whose Galois
representations are compatible with an incrementally relaxing
set of conditions.  As the conditions loosen, the set of forms
must grow to keep pace with the increasing number of lifts.
The increase in the number of lifts can be estimated from above
by a local cohomological calculation.  A sufficient supply of modular
forms is then furnished by the theory of congruences between
normalized
eigenforms of differing level.  This latter theory may be
viewed as a vast
generalization of what went on in
the author's article~\cite\ICM.

At the risk of distorting slightly the theory, I 
will 
index the shifting set of local conditions by 
a finite set of prime numbers~$\Sigma$ which
contains the set $\Sigma_0$
of primes at which
$\rho_{A,\ell}$ is ramified.
This latter set may be obtained concretely as
the set of primes which divide the discriminant of the
number field obtaining by adjoining to~$\Q$ the coordinates
of all points in~$A[\ell]$.  It is not hard to see that
$\Sigma_0$ contains~$\ell$
and is contained in the union of~$\{\ell\}$ and the set of
primes at which $A$ has bad reduction. Indeed,
this union is the
set of primes $\Sigma^\ast$ 
at which~$\rho_{A,\ell^\infty}$ is ramified,
according to
the well-known criterion of N\'eron-Ogg-Shafarevich~\cite\ST.
In this perspective, the minimal set of conditions
corresponds to the choice $\Sigma=\Sigma_0$, while a set
of conditions broad enough to include~$\rho_{A,\ell^\infty}$
is obtained by choosing $\Sigma=\Sigma^\ast$.

As promised, we shall now give the flavor of the
conditions that one imposes on the lifts~$\tilde\rho$
of type~$\Sigma$.
Firstly, we demand that $\tilde\rho$ be unramified outside~$\Sigma$.
Secondly, we ask that $\tilde\rho$ 
have the same qualitative
behavior at each prime $p\in\Sigma$
as the representation~$\rho_{A,\ell}$.
(The imposition of this condition can be traced
back to
Gouv\^ea's article~\cite\GOUCONTROL.)  Since $\ell$ lies
in~$\Sigma_0$,
a condition is imposed on~$\tilde\rho$
locally at the prime~$\ell$\emrule this condition requires
that $\tilde\rho$
be ``ordinary'' if $A$ has ordinary or multiplicative
reduction at~$\ell$ and that $\tilde\rho$ be ``flat''
if $A$ has supersingular reduction at~$\ell$, 
cf.~\cite\RUBINSILVER.  
For convenience, a
final (global) condition may be imposed
on~$\tilde\rho$, to the effect that
the determinant of~$\tilde\rho$
be the composite of the $\ell$-adic cyclotomic character
$$\tilde\chi_\ell\:\GalQ\to\Z_\ell^\ast$$
and the structural
map $\Z_\ell^\ast\to A^\ast$.  This supplementary
condition has the effect of allowing one
to work with the spaces $S(N)$ rather than
with spaces of modular forms on groups of the form~$\Gamma_1(N)$.

To prove that all lifts $\tilde\rho$ 
of type~$\Sigma$ are modular,
one first
passes to equivalence classes with respect to the
relation in which two representations
with values in~$\GL B$ are equivalent
whenever they are conjugate
by an element of~$\GL B$ which maps to the identity matrix
in~$\GL{\F_\ell}$.
The equivalence classes of lifts are called {\it deformations\/}\emrule
the terminology is borrowed from algebraic geometry, where
deformation theory has been developed extensively.
The idea of introducing deformation theory into the subject
of Galois representations is due to Mazur~\cite\MAZUR.
As F.~Gouv\^ea reports in his recent survey~\cite\GOUSURVEY,
the deformation viewpoint has evolved considerably since
\cite\MAZUR\ first appeared.
In particular, a 1993 article by R.~Ramakrishna~\cite\RAMA\
introduces foundational tools which
Wiles requires in~\cite\WILES; Ramakrishna shows
that deformation theory can be applied to study 
the family of lifts which
satisfy the local
conditions to which we have been alluding.

More precisely, the work of
Mazur and Ramakrishna
proves that there is a
{\it universal\/} deformation of type~$\Sigma$.
This is a lift
$$\rho_\Sigma\:\GalQ\to\GL{R_\Sigma}$$ which is
characterized by the property that for each lift
$$\tilde\rho\:\GalQ\to\GL B$$
of type~$\Sigma$
there is a unique homomorphism $\omega\:R_\Sigma\to B$
of local $\Z_\ell$-algebras
so that the deformation defined by~$\tilde\rho$ agrees with the
one obtained from~$\rho_\Sigma$ and~$\omega$.
A common
initial impression is that $R_\Sigma$ is
a relatively mysterious
object whose existence stems from an abstract representability
theorem\emrule one comes to grips with it only gradually.

Wiles seeks to compare $R_\Sigma$ with a concrete ring $\T_\Sigma$,
which he defines directly as a completion of a classical ring of
Hecke operators.  A theorem of Carayol constructs a Galois
representation
$$\rho'_\Sigma\:\GalQ\to\GL{\T_\Sigma}$$
with the property that the trace of
$\rho'_\Sigma(\Frob_p)$ is the Hecke operator $T_p$, for all
but finitely many primes~$p$.  One thinks of $\rho'_\Sigma$
as the
universal {\it modular\/} deformation
of type~$\Sigma$.
The problem is then
to prove that $\rho_\Sigma$ and~$\rho'_\Sigma$
coincide, so that all deformations of type~$\Sigma$
are modular.  
To come to grips with this problem,
Wiles begins with
the canonical
homomorphism $$\phi_\Sigma\: R_\Sigma\to\T_\Sigma$$ which
results from the universality of~$\rho_\Sigma$. 
It is relatively easy to show that $\phi_\Sigma$ is
surjective; the coincidence of
$\rho_\Sigma$ and~$\rho'_\Sigma$
means
that $\phi_\Sigma$ is an
{\it isomorphism}, cf.~\cite{\RUBINSILVER, \S4.2}.

My reaction to Wiles's 1993 announcement was
astonishment that one could prove
the modularity of Galois
representations by working directly with~$\phi_\Sigma$.
Despite conjectures by
Mazur~\cite{\MAZURTILOUINE, p.~85} and
Gouv\^ea~\cite{\GOUCONTROL, p.~108} to the effect
that $\phi_\Sigma$ is an isomorphism in the ordinary case,
and Theorem~2 of~\cite\FLACH,
I was not prepared for the revelation that
$\phi_\Sigma$
could be studied fruitfully.

\section Gorenstein and complete intersection conditions\endhead
The conjecture that $\phi_\Sigma$ is an isomorphism
is proved
in~\cite{\TW, \WILES}.
The proof
relies on standard notions of commutative
algebra which are discussed in~\cite\MATSUMURA\ and in~\cite\HWL.
The argument of Wiles and Taylor-Wiles
proceeds from the definition of~$\T_\Sigma$ as
a completion of the ring
generated by the Hecke operators $T_n$
acting on a specific space of classical cusp forms.
In particular, $\T_\Sigma$ is
free of finite rank over~$\Z_\ell$.
It has been known for some time that Hecke
rings such as~$\T_\Sigma$
look special from the vantage point of
commutative ring theory.  For example, $\T_\Sigma$ tends to
be Gorenstein, which means that the $\T_\Sigma$-module
$\Hom(\T_\Sigma,\Z_\ell)$ is free of rank~1
over~$\T_\Sigma$.  (The Gorenstein property was first noted in a
special case by Mazur \cite{\EIS, Ch.~II, \S15},
and then established in ever-widening generality by
others, including the author.)  Chapter~2
of~\cite\WILES\ includes a section which
summarizes and improves on the known Gorenstein assertions;
it proves, in particular,
that $\T_\Sigma$ is Gorenstein.

Using commutative algebra techniques,
\cite\WILES\ presents a number of conditions each of which is
sufficient to show
that $\phi_\Sigma$ is an isomorphism.
As stated in the Introduction above,
Wiles became aware
while studying the problem
that
$\phi_\Sigma$ is an isomorphism whenever the ring $\T_\Sigma$
is a complete intersection ring.
(This implication is proved from another point of view
in~\cite\MAZURNOTES.)
Later,
judging that it would be difficult to prove directly 
that $\T_\Sigma$ is a complete intersection ring, Wiles
focused on a numerical inequality~\cite{\RUBINSILVER, Th.~5.2}
and showed that
$\phi_\Sigma$ is an isomorphism whenever it is satisfied.

At the time of his 1993 Cambridge lectures,
Wiles believed that he had proved the numerical
inequality through the construction of
a ``geometric Euler system,'' thereby
generalizing work of~M.~Flach~\cite\FLACH.
Later analysis showed that the
construction envisaged by Wiles was flawed.
Interestingly,
it is not yet clear whether it can be modified so
as to yield an Euler system with the desired properties.
At the minimum, one
feels that this avenue of
inquiry is likely to remain extremely active.
In particular, Mazur's course notes~\cite\MAZURNOTES\
extract new
information from Flach's original
construction.

At the time of this revision,
the arguments of
\cite\TW\ and~\cite\WILES\
are being studied and internalized by
the mathematical community. 
Followup work is already beginning to appear.
For example, we mentioned above that
the main theorem of~\cite\WILES\ has been
strengthened by a recent manuscript 
of~F.~Diamond~\cite\DIAMONDNEW.  Also, the 
arguments given in
\cite\TW\ and~\cite\WILES\
have been shortened somewhat by
simplifications due
to~G.~Faltings; these
are explained in an appendix to~\cite\TW.

\section Toward the \Weil\ conjecture\endhead
We conclude with a short survey of results which have been
obtained by Wiles~\cite\WILES, Taylor-Wiles~\cite\TW, and
Diamond~\cite\DIAMONDNEW.
Two theorems which have the flavor
$$ \rho_{A,\ell} \hbox{ modular } \Longrightarrow
\rho_{A,\ell^\infty} \hbox{ modular }$$
are presented
in~\cite\WILES.
In each of the theorems, the prime $\ell$ is taken to be~odd,
and the representation $\rho_{A,\ell}$ is required to
be irreducible.  

One of the theorems (namely, \cite{\WILES, Th.~4.8})
has no application to Fermat's
Last Theorem, since it 
does not apply to semistable elliptic
curves.  This theorem does apply to
in situations where 
$\rho_{A,\ell}$
is isomorphic to the representation obtained
from the $\ell$-division
points of a
complex multiplication elliptic curve $A'$ over~$\Q$.
In these cases,
$A$ need not have complex multiplication itself; it is merely
``linked mod~$\ell$'' to a CM~curve.  For a discussion of
this theorem, including an explicit determination of
the family of curves
which can be linked to a fixed curve $A'$,
see~\cite\RSNOTE.

The theorem in~\cite\WILES\
which applies to Fermat's Last Theorem is the
one whose proof depends on the new work of Taylor-Wiles~\cite\TW:
\proclaim{Theorem 1}
Suppose that $A$ is a semistable elliptic curve over~$\Q$.
Let $\ell$ be an odd prime.
Assume that the representation $\rho_{A,\ell}$ is both irreducible
and modular.  
Then 
$A$ is a modular elliptic curve.
\endproclaim
\noindent
Fermat's Last Theorem may then be proved by combining
the author's theorem~\cite\FERMAT\
with the following result, which may be viewed as a
highly non-obvious corollary of
Theorem~1.
\proclaim{Theorem 2}
Let $A$ be a semistable elliptic curve over~$\Q$.
Then $A$ is a modular elliptic curve.
\endproclaim
\noindent
Wiles deduces Theorem 2 from Theorem 1 by
an ingenious argument, which we will now describe.
(The argument, which
is presented in~\cite{\WILES, Ch.~5},
has been sketched in~\cite\RUBINSILVER.)
Let $A$ be a semistable elliptic curve, and consider the
representation
$\rho_{A,3}$.
If this representation happens to be irreducible, then 
it is also modular, by the
the
results of Langlands and Tunnell which were cited above.
Thus Theorem~1,
proves that $A$ is modular.

What happens if $\rho_{A,3}$ is reducible?
In this case, we
examine~$\rho_{A,5}$.  
If this latter representation is
reducible as well, then Wiles shows
directly that $A$ is modular.
Hence we can, and do,
suppose that $\rho_{A,5}$ is irreducible.
Wiles shows then that one can find
a second semistable elliptic
curve $A'$
whose mod~5 representation is isomorphic to that of~$A$ and
whose mod~3 representation is irreducible
(cf.~\cite{\RUBINSILVER,
Appendix~B}).
Two applications of Theorem~1
then suffice to show that $A$ is modular.
Indeed, applying the theorem to~$A'$ with $\ell=3$, we find
that $A'$ is modular.  In particular, the irreducible
representation $\rho_{A',5}$ is modular.  Since this
representation coincides with~$\rho_{A,5}$, we may apply
the theorem to~$A$ with $\ell=5$ to conclude that $A$
is modular, as desired.

A preprint
of~F.~Diamond~\cite\DIAMONDNEW\ generalizes 
Theorem~1
to 
the case where $A$ is an elliptic curve over~$\Q$ whose
conductor is not divisible by~$\ell^2$
(i.e, one which is semistable at~$\ell$):

\proclaim{Theorem 3}
Let $\ell$ be an odd prime number.
Suppose that
$A$ is an elliptic curve over~$\Q$ 
which is semistable at~$\ell$.
Assume that the representation $\rho_{A,\ell}$ is both irreducible
and modular.  
Further, if $\ell=3$, assume
that
$\rho_{A,3}\left(\Gal(\Qbar/\Q\left(\sqrt{-3}\right))\right)$
is non-abelian.  Then 
$A$ is a modular elliptic curve.
\endproclaim
\noindent
The condition concerning 
$\Gal\left(\Qbar/\Q\left(\sqrt{-3}\right)\right)$
occurs already in Wiles's work.  However, in the situation
of Theorem~1, Wiles proves that
$\rho_{A,3}\left(\Gal\left(\Qbar/\Q\left(\sqrt{-3}\right)\right)\right)$
is irreducible and non-abelian when~$\ell=3$; in other
words, the assumption relative 
to~$\Gal\left(\Qbar/\Q\left(\sqrt{-3}\right)\right)$ has
been omitted from Theorem~1 because it may be proved
unconditionally.

Using a variant of the Wiles
argument we have just sketched, Diamond
deduces the following generalization of
Theorem~2.
\proclaim{Theorem 4}
Suppose that $A$ is an elliptic curve over~$\Q$ which is
semistable both at~$3$ and at~$5$.  Then $A$ is modular.
\endproclaim


\Refs
\catcode`\?=\active
\def?{.\hskip 0.1667em\relax}

\ref\no\ASHGROSS\by A. Ash and R. Gross
\paper From Gauss
to Langlands: a context for Wiles's achievement
\paperinfo in preparation
\endref

\ref\no\AL\by A?O?L. Atkin and J. Lehner
\paper
Hecke operators on $\go m$
\jour Math. Annalen
\vol 185\yr1970\pages134--160
\endref

\ref\no\GP\by C. Batut, D. Bernardi, H. Cohen and M. Olivier
\paper GP/PARI
\paperinfo Available by anonymous ftp
from {\tt megrez.math.u-bordeaux.fr} or
{\tt math.ucla.edu}, in the directory
{\tt /pub/pari}
\endref

\ref\no\ANTWERPIV\by B?J. Birch and W. Kuyk, eds.
\book Modular functions of one variable IV
\bookinfo Lecture Notes in Math., vol.~476
\publ Springer-Verlag
\publaddr Berlin and New York
\yr1975\endref

\ref\no\BOSMA\by W. Bosma and H?W. Lenstra, Jr.
\paper Complete systems of two addition laws
for elliptic curves
\toappear
\endref

\ref\no\BOSTONSURVEY\by N. Boston
\paper A Taylor-made plug for Wiles' proof
\jour College Math. J.
\toappear
\endref

\ref\no\BOSTONGRANVILLE\by N. Boston and A. Granville
\paper Review of \cite\MvS
\jour American Math. Monthly
\toappear
\endref

\ref\no\BUZ\by K?M. Buzzard
\paper The levels of modular representations
\paperinfo Cambridge University thesis, 1995
\endref

\ref\no\CAREARLY \by H. Carayol
\paper Sur les repr\'esentations
  $\ell$-adiques associ\'ees aux formes modulaires de Hilbert
\jour  Ann. scient. \'Ec. Norm. Sup., $4^{\roman e}$ s\'erie
\vol 19 \yr1986 \pages409--468\endref 

\ref\no\COATESIP\by J?H. Coates and S?T. Yau, eds.
\book Elliptic curves and modular forms
\bookinfo proceedings of a conference held in Hong Kong,
December 18--21, 1993
\publ International Press
\publaddr Cambridge, MA and Hong Kong
\toappear
\endref

\ref\no\CONNELL\by I. Connell
\paper Apecs (arithmetic of plane elliptic curves)\emrule a
program written in Maple
\paperinfo Available by anonymous ftp from {\tt math.mcgill.ca}
in the directory {\tt /pub/apecs}
\endref

\ref\no\CORNELLSILVERMAN\by G. Cornell and J. Silverman, eds.
\book Arithmetic Geometry
\publ Springer-Verlag
\publaddr Berlin and New York
\yr1986\endref

\ref\no\COX\by D. Cox
\paper Introduction to Fermat's Last Theorem
\jour American Math. Monthly
\vol101\yr1994\pages3--14
\endref

\ref\no\DARMON\by H. Darmon
\paper 
The Shimura-Taniyama conjecture (d'apr\`es Wiles)
\jour Russian Math Surveys
\toappear
\endref

\ref\no\NEWDARMON\bysame
\paper Serre's conjectures
\paperinfo in~\cite\MURTYASED
\endref

\ref\no\DELIGNESERRE\by P. Deligne and J.-P. Serre
   \paper Formes modulaires de poids~$1$
\jour
   Ann. scient. \'Ec. Norm. Sup., $4^{\roman e}$ s\'erie
     \vol 7\yr1974\pages 507--530
\endref

\ref\no\DIAMONDREF\by F. Diamond
\paper The refined conjecture of Serre
\paperinfo in~\cite\COATESIP
\endref

\ref\no\DIAMONDNEW\bysame
\paper On deformation rings and Hecke rings
\toappear
\endref

\ref\no\FALT\by G. Faltings
\paper Endlichkeitss\"atze f\"ur abelsche
    Variet\"aten \"uber Zahl\-k\"orpern
\jour Invent. Math.\vol73\yr1983\pages349--366\endref

\ref\no\FLACH\by
M. Flach
\paper
A finiteness theorem for the symmetric square of an elliptic
curve
\jour Invent. Math.\vol 109\yr1992\pages307--327
\endref

\ref\no\FONT\by J.-M. Fontaine
\paper Il n'y a pas de vari\'et\'e ab\'elienne sur $\Z$
\jour Invent. Math.\vol81\yr1985\pages515--538
\endref

\ref\no\FONTAINEMAZUR
\by J.-M. Fontaine and B. Mazur
\paper Geometric Galois representations
\paperinfo in~\cite\COATESIP
\endref

\ref\no\FREYONE\by G. Frey
\paper Links between stable elliptic curves and certain
  diophantine equations
\jour Annales Universitatis Saraviensis
\vol1\yr1986\pages1--40
\endref

\ref\no\FREYTWO\bysame
\paper Links between elliptic curves and solutions
of $A-B=C$\jour Journal of the Indian Math. Soc.\vol51
\yr1987\pages117-145\endref

\ref\no\GELBART\by S. Gelbart
\paper Automorphic forms and Artin's conjecture
\jour \jour Lecture Notes
in Math
\vol 627\yr1977\pages241--276\endref

\ref\no\LABESSE\by P. G\'erardin and J?P. Labesse
\paper The solution to a base change problem for
$\GLTWO$ (following Langlands, Saito, Shintani)
\jour
Proceedings of Symposia in Pure Mathematics
\vol33 (2)
\yr1979
\pages115--133
\endref

\ref\no\GOUCONTROL\by F?Q. Gouv\^ea
\paper Deforming Galois representations: controlling
the conductor
\jour Journal of Number Theory
\vol34\yr1990\pages95--113
\endref

\ref\no\GOUVEA\bysame
\paper ``A marvelous proof''\jour American Math. Monthly
\vol101\yr1994\pages203--222
\endref

\ref\no\GOUSURVEY\bysame
\paper Deforming Galois representations: a survey
\paperinfo in~\cite\MURTYASED
\endref

\ref\no\HAYESRIBET\by B. Hayes and K?A. Ribet
\paper Fermat's Last Theorem and modern arithmetic
\jour American Scientist
\vol82\yr1994
\pages144--156
\endref

\ref\no\HEARSTRIBET\by W. R. Hearst III and K?A. Ribet
\paper
Review of
``Rational points on elliptic curves'' by Joseph H. Silverman and 
 John T. Tate
\jour
Bulletin of the AMS
\vol30
\yr1994
\pages248--252
\endref

\ref\no\HINDRY\by M. Hindry\paper
``a, b, c'', conducteur, discriminant
\jour
Publications math\'ematiques de l'Univer\-sit\'e
Pierre et Marie Curie, Probl\`emes diophantiens
\yr1986--87
\endref

\ref\no\AXJ\by A. Jackson
\paper Update on proof of Fermat's Last Theorem
\jour
Notices of the AMS\vol41\yr1994
\pages185--186
\endref

\ref\no\ANOTHERAXJ\bysame
\paper Another step toward Fermat
\jour
Notices of the AMS\vol42\yr1995\pages48
\endref

\ref\no\KNAPP\by A?W.~Knapp
\book Elliptic curves
\bookinfo Math. Notes, vol. 40
\publ Princeton Univ. Press
\publaddr Princeton, NJ\yr1992
\endref

\ref\no\LANGBOOK\by S. Lang
\book Introduction to modular forms
\publ Springer-Verlag
\publaddr Berlin and New York
\yr1976\endref

\ref\no\LANGABELIANVAR\bysame
\book Abelian varieties
\publ Springer-Verlag
\publaddr Berlin and New York
\yr1983\endref

\ref\no\LANGFILE\bysame
\paper
The \Weil\ file
\paperinfo Available directly from S.~Lang, Yale Math.\ Department
\endref

\ref\no\LANG\bysame
\paper
Old and new conjectured diophantine inequalities
\jour Bull. AMS
\vol23
\yr1990
\pages37--75\endref

\ref\no\LANGLANDS\by R?P. Langlands
\book Base change for $\GLTWO$
\bookinfo Annals of Math. Studies, vol. 96
\publ Princeton University Press
\publaddr Princeton
\yr1980\endref

\ref\no\HWL\by H?W. Lenstra, Jr.
\paper Complete intersections and Gorenstein rings
\paperinfo in~\cite\COATESIP
\endref

\ref\no\LI\by W.-C? W. Li
\paper Newforms and functional equations
\jour Math. Annalen
\vol212\yr1975\pages285--315
\endref

\ref\no\MATSUMURA\by H. Matsumura
\book Commutative ring theory
\publ Cambridge University Press
\publaddr Cambridge (U.K.)
\yr1986
\endref

\ref\no\SILVERMAN\by P?A. van Mulbregt and J?H. Silverman
\paper Elliptic curve calculator
\paperinfo Available by anonymous ftp from
{\tt gauss.math.brown.edu}
in the directory {\tt /dist/EllipticCurve}
\endref

\ref\no\EIS\by B. Mazur
\paper Modular curves and the Eisenstein ideal
\jour Publ. Math. IHES
\vol47
\yr1977
\pages33--186
\endref

\ref\no\PRIMEDEGREE\bysame
\paper Rational isogenies of prime degree
\jour
Invent. Math.\vol44\yr1978\pages129--162\endref

\ref\no\MAZUR\bysame
\paper
Deforming Galois representations
\inbook
Galois groups over~$\Q$
\bookinfo MSRI Publications, vol. 16
\publ Springer-Verlag\publaddr Berlin and New York \yr1989
\pages385--437\endref

\ref\no\GAD\bysame
\paper Number theory as gadfly
\jour Am. Math. Monthly \vol98\yr1991\pages593--610
\endref

\ref\no\MAZURNOTES\bysame
\book Very rough course notes for Math 257y, parts I--III
\bookinfo to appear as
``Galois deformations and Hecke curves''
\endref

\ref\no\MAZURNUMBER\bysame
\paper Questions about number
\inbook New directions in mathematics
\toappear
\endref

\ref\no\MAZURTILOUINE\by B. Mazur and J. Tilouine
\paper Repr\'esentations galoisiennes, diff\'erentielles
de K\"ahler et \<<con\-jectures principales\>>
\jour Publ. Math. IHES
\vol71\yr1990
\pages9--103
\endref

\ref\no\MIYAKEOLD\by T. Miyake
\paper On automorphic forms on $\varGLTWO$ and
Hecke operators
\jour Annals of Math.\vol94\yr1971\pages174--189
\endref

\ref\no\MIYAKE\bysame
\book Modular forms
\publ Springer-Verlag
\publaddr Berlin and New York
\yr1989\endref

\ref\no\MURTY\by V?K. Murty
\book
Introduction to Abelian varieties
\bookinfo CRM monograph series, vol.~3
\publ American Mathematical Society
\publaddr Providence
\yr1993
\endref

\ref\no\MURTYASED\bysame, ed.
\book Elliptic curves, galois representations
and modular forms
\bookinfo CMS Conference Proceedings
\publ American Mathematical Society
\publaddr Providence
\toappear
\endref

\ref\no\OESTERLE\by J. Oesterl\'e
\paper Nouvelles approches du ``th\'eor\`eme'' de Fermat
\jour Ast\'erisque
\yr 1988
\vol 161/162
\pages165--186
\endref

\ref\no\OGG\by A. Ogg
\paper
Elliptic curves and wild ramification
\jour
American J. of Math.
\yr1967
\vol89
\pages1--21
\endref

\ref\no\PRASAD\by D. Prasad
\paper Ribet's Theorem: Shimura-Taniyama-Weil implies Fermat
\paperinfo in~\cite\MURTYASED
\endref

\ref\no\RAMA\by R. Ramakrishna
\paper On a variation of Mazur's deformation functor
\jour Compositio Math.
\vol87\yr1993\pages269--286
\endref

\ref\no\RIBETSURVEY\by K?A. Ribet\paper
The $\ell$-adic representations attached to
an eigenform with Nebentypus: a survey
\jour Lecture Notes
in Math\vol 601\yr1977\pages17--52\endref

\ref\no\ICM\bysame
\paper  Congruence relations between modular forms
\inbook Proc. Int. Cong. of Mathematicians 1983
\pages503--514\endref

\ref\no\FERMAT\bysame
\paper
On modular representations of $\GalQ$
arising from modular forms\jour
Invent. Math.\vol100\yr1990\pages431--476\endref

\ref\no\TOULOUSE\bysame
\paper From the \Weil\ conjecture to Fermat's Last
Theorem
\jour Annales de la Facult\'e des Sciences de l'Universit\'e
   de Toulouse
\vol 11
\yr1990
\pages116--139
\endref

\ref\no\KOREA\bysame\paper
Abelian varieties over $\Q$ and modular forms
\inbook
1992 Proceedings
of KAIST Mathematics Workshop
\publaddr Taejon
\publ Korea Advanced Institute of Science and Technology
\yr 1992\pages53--79
\endref

\ref\no\NOTICES\bysame
\paper Wiles proves Taniyama's conjecture; Fermat's
Last Theorem follows\jour Notices of the AMS\vol40\yr1993
\pages575--576
\endref

\ref\no\VIDEO\bysame
\book Modular elliptic curves and Fermat's last theorem:
a lecture presented at 
George Washington University, Washington, DC,
August 1993
\bookinfo Selected Lectures in Mathematics
\publ American Mathematical Society
\publaddr Providence
\yr1993
\endref

\ref\no\MOTIVES\bysame
\paper
Report on mod~$\ell$ representations of
$\GalQ$
\jour
Proceedings of Symposia in Pure Mathematics
\vol55 (2)
\yr1994
\pages639--676
\endref

\ref\no\RUBINSILVER\by K. Rubin and A. Silverberg
\paper A report on Wiles' Cambridge Lectures
\jour Bulletin of the AMS
\yr1994
\vol31
\pages15--38
\endref

\ref\no\RSNOTE\bysame
\paper
Families of elliptic curves with constant mod~$p$
representations
\paperinfo in~\cite\COATESIP
\endref

\ref\no\SERRETAU\by J-P. Serre
\paper
Une interpr\'etation des congruences relatives \`a
la fonction $\tau$ de Ramanujan
\jour
S\'eminaire Delange-Pisot-Poitou
\yr 1967--68, ${\roman n}^{\roman o}$ 14
\endref

\ref\no\SERREOLD\bysame
\paper
Facteurs locaux des fonctions z\^eta des
vari\'et\'es alg\'ebriques (d\'efinitions
et conjectures)
\jour
S\'eminaire Delange-Pisot-Poitou
\yr 1969--70, ${\roman n}^{\roman o}$ 19
\endref

\ref\no\SERREIM\bysame\paper Propri\'et\'es galoisiennes
 des points d'ordre fini des courbes elliptiques
\jour
 Invent. Math.\vol 15\yr1972\pages259--331\endref

\ref\no\SERRECOURS\bysame
\book A course in arithmetic
\bookinfo Graduate Texts in Math., vol. 7
\publ Springer-Verlag \publaddr New York, Heidelberg and Berlin
\yr1973\endref

\ref\no\ARCATA\bysame
\paper Lettre \`a J.-F. Mestre (13 ao\^ut 1985)\jour
 Contemporary Mathematics \vol 67\yr1987\pages263--268\endref

\ref\no\DUKE\bysame
\paper Sur les repr\'esentations modulaires de degr\'e~2
de $\GalQ$\jour Duke Math. J.\vol 54\yr1987\pages179--230\endref

\ref\no\GROUPESALG\bysame
\book Algebraic groups and class fields
\bookinfo Graduate Texts in Math., vol. 117
\publ Springer-Verlag \publaddr New York, Heidelberg and Berlin
\yr1988\endref


\ref\no\ST\by J-P. Serre and J. Tate
\paper Good reduction of abelian varieties
\jour Ann. of Math. \vol 88\yr1968\pages492--517\endref

\ref\no\SHIMURAOLD\by G. Shimura
\paper An $\ell$-adic method in the theory of automorphic forms
\endref

\ref\no\SHCRELLE\bysame
\paper A reciprocity law in non-solvable extensions
\jour Journal f\"ur die reine und angewandte Mathematik
\vol221\yr1966\pages209--220
\endref

\ref\no\SHIMURABOOK\bysame
\book Introduction to the arithmetic theory of automorphic functions
\publ Princeton University Press
\publaddr Princeton \yr1971\endref

\ref\no\SHIMURAPINK\bysame
\paper
On elliptic curves with complex
multiplication as factors of the Jacobians of modular
function fields\jour Nagoya Math. J.\vol43\yr1971
\pages199--208\endref

\ref\no\SHIMURAHECKE\bysame
\paper Class fields over real quadratic fields and
Hecke operators
\jour Ann. of Math.\vol95\yr1972\pages131--190\endref

\ref\no\SHIMURAFACTORS\bysame
\paper On the factors of the jacobian variety of a
  modular function field
\jour J. Math. Soc. Japan
\vol25
\yr1973
\pages523--544\endref

\ref\no\SILVERMANOLD\by J?H. Silverman
\book The arithmetic of elliptic curves
\bookinfo Graduate Texts in Math., vol.~106
\publ Springer-Verlag\publaddr Berlin and New York
\yr1986
\endref

\ref\no\SILVERMANNEW\bysame
\book Advanced topics in the arithmetic of elliptic curves
\bookinfo Graduate Texts in Math., vol.~151
\publ Springer-Verlag
\publaddr Berlin and New York
\yr1994
\endref

\ref\no\TATE\by J?T. Tate
\paper
Algorithm for determining the type of a singular
fiber in an elliptic pencil
\jour
Lecture Notes
in Math\vol 476\yr1975\pages33--52
\endref

\ref\no\TATELETTER\bysame
\paper The non-existence of certain Galois extension of~$\Q$
unramified outside~2
\jour Contemporary Mathematics
\vol174\yr1994\pages153--156\endref

\ref\no\TW\by R?L. Taylor and A. Wiles
\paper
Ring theoretic properties of certain Hecke algebras
\jour
Annals of Math.
\toappear
\endref

\ref\no\TUNNELL\by
J. Tunnell
\paper
Artin's conjecture for representations
of octahedral type
\jour Bull. AMS (new series)
\vol5
\yr1981
\pages173--175
\endref

\ref\no\MvS\by M. vos Savant
\book The world's most famous math problem:
the proof of Fermat's last theorem and other mathematical
mysteries
\publ St.~Martin's Press
\publaddr New York\yr1993
\endref

\ref\no\WEIL\by A. Weil
\paper \"Uber die Bestimmung Dirichletscher Reihen
  durch Funktionalgleichungen
\jour Math. Annalen\vol168\yr1967\pages165--172
\endref

\ref\no\WILES\by A. Wiles
\paper Modular elliptic curves and Fermat's Last Theorem
\jour Annals of Math.
\toappear
\endref

\endRefs
\enddocument